# ESTIMATION OF FRACTAL DIMENSION FOR A CLASS OF NON-GAUSSIAN STATIONARY PROCESSES AND FIELDS

BY GRACE CHAN AND ANDREW T. A. WOOD

*University of Iowa and University of Nottingham*

We present the asymptotic distribution theory for a class of increment-based estimators of the fractal dimension of a random field of the form $g\{X(t)\}$, where $g:\mathbf{R}\to\mathbf{R}$ is an unknown smooth function and $X(t)$ is a real-valued stationary Gaussian field on $\mathbf{R}^d$, $d=1$ or 2, whose covariance function obeys a power law at the origin. The relevant theoretical framework here is "fixed domain" (or "infill") asymptotics. Surprisingly, the limit theory in this non-Gaussian case is somewhat richer than in the Gaussian case (the latter is recovered when $g$ is affine), in part because estimators of the type considered may have an asymptotic variance which is random in the limit. Broadly, when $g$ is smooth and nonaffine, three types of limit distributions can arise, types (i), (ii) and (iii), say. Each type can be represented as a random integral. More specifically, type (i) can be represented as the integral of a certain random function with respect to Lebesgue measure; type (ii) can be represented as the integral of a second random function with respect to an independent Gaussian random measure; and type (iii) can be represented as a Wiener–Itô integral of order 2. Which type occurs depends on a combination of the following factors: the roughness of $X(t)$, whether $d=1$ or $d=2$ and the order of the increment which is used. Another notable feature of our results is that, even though the estimators we consider are based on a variogram, no moment conditions are required on the observed field $g\{X(t)\}$ for the limit theory to hold. The results of a numerical study are also presented.

## 1. Introduction.

1.1. *Background.* The problem of quantifying the roughness of a (continuous but rough) curve or surface, whose height is observed at discrete locations on a rectangular grid, arises in many areas of science and technology. A widely used approach is to model the curve or surface as a random field whose covariance function follows some form of power law behavior at the origin, and then to estimate a scale-invariant measure such as the fractal dimension to quantify roughness. A good entry point to recent statistical literature on this topic is the discussion paper by Davies and Hall (1999).

Frequently in this approach it has been assumed that the curve or surface is Gaussian. However, in a number of applications the Gaussian assumption may be open to doubt, and therefore the problem of estimating the fractal dimension in non-Gaussian settings is of interest. In this paper we study estimators of fractal dimension for a class of stationary non-Gaussian random fields, and we provide a detailed account of the asymptotic distribution theory for these estimators, as well as studying their numerical properties using simulation. As a preliminary, we provide a brief review of recent work in which the Gaussian assumption is made.

1.2. *The stationary Gaussian model.* For simplicity we focus mainly on the one-dimensional case ($d = 1$) in this Introduction. Let $\{X(t) : t \in \mathbf{R}\}$ denote a stationary Gaussian process with a covariance function $\gamma$ which obeys the following power law at the origin:

$$(1.1) \qquad \gamma(t) = \gamma(0) - c|t|^\alpha + O(|t|^{\alpha+\beta}) \qquad \text{as } |t| \to 0,$$

where $\alpha \in (0, 2]$, known as the fractal index, governs the roughness of the sample functions and is typically the parameter of greatest interest in roughness studies; the positive quantity $c$ is a (local) scale parameter known as the topothesy; and $\beta > 0$ governs the size of the remainder term in (1.1). There is a simple relationship, under model (1.1), between $\alpha$ and the fractal dimension $D$ of the graph of the sample function, given by $D = 2 - \alpha/2$; see, for example, Adler (1981). [This result generalizes to $D = d + 1 - \alpha/2$ when $X(t)$ is a stationary Gaussian random field on $\mathbf{R}^d$ with covariance function of the form (2.3).] Thus, the larger the value of $\alpha$, the smoother the sample function.

Suppose we observe a sample

$$(1.2) \qquad \mathcal{S}_n = \{X(i/n) : i = 0, 1, \ldots, n-1\}$$

of observations of $X(t)$ at equally spaced locations in the region $[0, 1]$. In this formulation of the problem the asymptotic regime as $n \to \infty$ is known as "fixed domain asymptotics" and is often appropriate when interest is focused on roughness at fine scales. See Stein (1999) for further discussion of this type of asymptotic regime.

A variety of estimators of $\alpha$, based on data $\mathcal{S}_n$, have been studied in recent years under model (1.1), with $X$ assumed to be stationary and Gaussian.

For example, Hall and Wood (1993) considered box-counting estimators; Jakeman and Jordan (1990) and Constantine and Hall (1994) discussed estimators based on the variogram; Feuerverger, Hall and Wood (1994) considered estimators based on counting upcrossings; and Chan, Hall and Poskitt (1995) considered estimators based on the periodogram. Kent and Wood (1997) considered two modifications to variogram-based estimators: the use of higher-order increments and the use of generalized least squares. The use of higher-order increments in a closely related context was investigated independently by Istas and Lang (1997).

In this paper the focus is on variogram-based estimators of fractal dimension, and we now discuss these in more detail. The theoretical variogram is given by $\nu(h) = E\{X(t+h) - X(t)\}^2$ and under (1.1) $\nu(h) \sim 2c|h|^\alpha$ as $|h| \to 0$. Given data of the form (1.2), we can estimate $\nu(u/n)$ by

$$\hat{\nu}_u = \sum_i [X\{(i+u)/n\} - X(i/n)]^2, \tag{1.3}$$

and we can estimate $\alpha$ using log–log regression based on the approximation relationship

$$\log \hat{\nu}_u \approx \text{const.} + \alpha \log u. \tag{1.4}$$

The simplest estimator of this type is the ordinary least squares estimator given by

$$\hat{\alpha} = \frac{\sum_{u=1}^m \log \hat{\nu}_u (\log u - m^{-1} \sum_{k=1}^m \log k)}{\sum_{u=1}^m (\log u - m^{-1} \sum_{k=1}^m \log k)^2}, \tag{1.5}$$

where $m$ stays fixed as $n \to \infty$.

It turns out that when $0 < \alpha < 3/2$ the estimator (1.5) has variance of order $n^{-1}$, but when $3/2 < \alpha < 2$ the variance of (1.5) is of order $n^{2\alpha-4}$. See Constantine and Hall (1994); related results were also obtained by Jakeman and Jordan (1990). More recently, it was noted by Kent and Wood (1997) and Istas and Lang (1997) that if one bases the variogram in (1.3) on second-order differences, that is,

$$\hat{\nu}_u = \sum_i [X\{(i+u)/n\} + X\{(i-u)/n\} - 2X(i/n)]^2, \tag{1.6}$$

then the variance of the resulting estimator (1.5) is of order $n^{-1}$ for all $\alpha \in (0, 2)$. Thus, there is motivation for considering variograms based on increments (i.e., differences) of higher order. In our terminology (1.3) is based on an increment of order 0, and (1.6) is based on an increment of order 1, and in either case (1.3) may be written in the form

$$\hat{\nu}_u = \sum_i \left\{ \sum_j a_j^u X\left(\frac{i+j}{n}\right) \right\}^2,$$

where the $a_j^u$ notation for increments is described in detail in Kent and Wood (1997) and Chan and Wood (2000) and is summarized in the Appendix.

More recently, the corresponding estimation problem for the two-dimensional case was considered in the discussion paper by Davies and Hall (1999). The analysis of data obtained by two-dimensional sampling is trickier than the one-dimensional case, in part because of the possibility that anisotropy is present. A helpful discussion of this issue is given by Davies and Hall (1999). A second approach to the analysis of two-dimensional surface data is considered by Chan and Wood (2000).

1.3. *Non-Gaussian data.* In this paper we assume that, instead of observing the stationary Gaussian process $X(t)$, we observe a stationary non-Gaussian process $g\{X(t)\}$, where $g$ is a smooth but unknown nonaffine function and $X$ is a stationary Gaussian process satisfying (1.1) as before. Two-dimensional sampling of random fields (corresponding to $d = 2$) is also considered. We address the following question: what are the asymptotic properties of the variogram-based estimators when we observe $g\{X(t)\}$ rather than $X(t)$? From a theoretical point of view, the quantity being estimated, namely the fractal dimension of the sample function, will be the same as that of the underlying Gaussian random field provided $g$ is nondegenerate and reasonably smooth; see Hall and Roy (1994) for the relevant results. However, it turns out that the asymptotic distribution theory for nonaffine $g$ is somewhat richer than in the Gaussian case (though the Gaussian case is recovered, of course, when $g$ is an affine transformation).

Chan and Wood (2000) obtained (correct) preliminary results concerning rates of convergence in the case when $g(x) = x^2$. However, no concrete results concerning the asymptotic distribution theory were given, and in fact the conjectures in Remarks 5.5 and 5.6 of that paper do not adequately describe the results given below. Also, it turns out that the case $g(x) = x^2$ is not fully representative of general smooth nondegenerate $g$, because some components of the limit distribution disappear or are constant in the former case.

Our main theoretical results are stated in Section 2. One general point which emerges is that the estimation of fractal dimension in this non-Gaussian setting is, in a sense, more difficult than in the Gaussian case; see, in particular, point 2 in Section 2.3. However, on the numerical side a fairly extensive simulation study reported in Section 3 suggests that, in practice, the deterioration in the non-Gaussian setting is fairly mild, provided that $g$ is not drastically nonaffine in the relevant domain. The main theorems are proved in Section 5. These proofs make use of several lemmas which are proved in Section 4.

**2. Main results.** The principal results in the paper are presented in Theorems A and B. Theorem A covers those cases in which $4 + 4p - 2\alpha > d$ and Theorem B covers those cases in which $4 + 4p - 2\alpha < d$, where $p \geq 0$ is the order of the increment used (see the Appendix), $d = 1, 2$ is the dimension of the parameter set of the underlying random field and $\alpha \in (0, 2)$ is the fractal index.

Throughout this paper we make frequent use of the notation for increments and multi-indices given in the Appendix.

2.1. *Preliminaries.* Let $\{X(t) : \in \mathbf{R}^d\}$ denote a real-valued stationary Gaussian process $(d = 1)$ or field $(d = 2)$ with covariance function given by

$$\gamma(t) = \text{cov}\{X(s), X(s+t)\}.$$

Let $g : \mathbf{R} \to \mathbf{R}$ be an unknown function. It is assumed that we observe $g\{X(t)\}$ rather than $X(t)$.

Define the index set

(2.1) $$\mathcal{I}_n = \{0 \leq j < n_0\}.$$

When $d = 1$, $j$ and $n_0$ are integers and we take $n = n_0$; and when $d = 2$, $j$ and $n_0$ are multi-indices in $\mathbf{Z}^2$ (see the Appendix) and $n_0 = n_0(n)$ is a sequence such that $n$ is the product of the elements of $n_0$, that is, $n = n_0[1]n_0[2]$. Thus, $n$ is the number of elements in $\mathcal{I}_n$. The dataset we actually sample is given by

$$\mathcal{S}_n = \{g_i \equiv g\{X(i/n_0)\} : -mJ \leq i < n_0 + mJ\},$$

where division of multi-indices $j, k \in \mathbf{Z}^d$ is defined by $j/k = (j[1]/k[1], \ldots, j[d]/k[d])$, assuming $k[l] \neq 0$ for each $l$. In the definition of $\mathcal{S}_n$, $J$ is a multi-index which depends on the increment that we use, and $m$, an integer, is the number of dilations of the increment that we consider [see the Appendix and also Kent and Wood (1997)]. It is assumed throughout that $m$ stays fixed as $n \to \infty$; justification for keeping $m$ fixed is given by Constantine and Hall (1994). Note that the sampling regime indicated above corresponds to "fixed domain" asymptotics as $n \to \infty$.

Consider the following conditions on $\gamma$ and $g$.

$(\mathcal{A}1)_q^{(d)}$ For some $\alpha \in (0, 2)$ and $\beta > 0$, and for each nonnegative multi-index $r$ with $|r| = q$,

(2.2) $$\gamma^{(r)}(t) = -\frac{\partial^q}{\partial t^r}\{\|t\|^\alpha M(t/\|t\|)\} + O(\|t\|^{\alpha+\beta-q})$$

$$\text{as } \|t\| \to 0,$$

where $\|t\| = (t^T t)^{1/2}$ is the usual Euclidean norm on $\mathbf{R}^d$, and, for a nonnegative multi-index $r = (r[1], \ldots, r[d])$, $\gamma^{(r)}(t) = \partial^{|r|}\gamma(t)/\partial t^r =$ $\partial^{|r|}\gamma(t)/\partial t_1^{r[1]}\cdots\partial t_d^{r[d]}$, where $|r| = \sum_{l=1}^{d} r[l]$. In (2.2) $M(\cdot)$ is assumed to be a positive constant when $d=1$; when $d=2$ $M(\cdot)$ is a function on the unit circle in $\mathbf{R}^2$, all of whose partial derivatives derivatives of order $q+1$ are assumed continuous.

($\mathcal{A}2$) The seventh derivative of $g$, $g^{(7)}$, is continuous on $\mathbf{R}$.

($\mathcal{A}3$) The set $\{x : g^{(1)}(x) = 0\}$ has Lebesgue measure 0, where $g^{(1)}$ is the derivative of $g$.

($\mathcal{A}4$) ($d=2$ only.) As $n \to \infty$ $n_0[1]/n_0[2]$ stays bounded away from 0 and $\infty$.

It is easily shown that, if $(\mathcal{A}1)_q^{(d)}$ holds, then $(\mathcal{A}1)_r^{(d)}$ holds for each $1 \le r < q$. Moreover, $(\mathcal{A}1)_q^{(d)}$ implies that

$$(2.3) \qquad \gamma(t) = \gamma(0) - \|t\|^\alpha M(t/\|t\|) + O(\|t\|^{\alpha+\beta}),$$

where (2.3) reduces to (1.1) when $d=1$.

It is possible to weaken assumption ($\mathcal{A}2$) to some extent, but only at a considerable cost in technical detail in the proofs. Assumption ($\mathcal{A}3$) is a mild nondegeneracy condition which seems essential if our theorems are to hold. Condition ($\mathcal{A}4$) is needed when $d=2$ to ensure that the sampling set does not become too "thin."

Putting $g_i = g\{X(i/n_0)\}$, we define

$$(2.4) \qquad \bar{Z}_u = n^{-1} \sum_{i \in \mathcal{I}_n} \left\{ \sum_j a_j^u g_{i+j} \right\}^2, \qquad \mu_u = n^{\alpha/d} E\left\{ \sum_j a_j^u X\left(\frac{i+j}{n_0}\right) \right\}^2$$

and

$$(2.5) \qquad G_r = \int_{t \in [0,1]^d} [g^{(1)}\{X(t)\}]^{2r}\, dt, \qquad r=1,2,$$

where $d=1$ or 2. The notation for increments $\mathbf{a} = \{a_j\}$ used in (2.4) is explained in the Appendix. Let $L_u$, $u=1,\ldots,m$, be real numbers which satisfy

$$\sum_{u=1}^m L_u = 0, \qquad \sum_{u=1}^m L_u \log u = 1.$$

Various choices for the $L_u$ are discussed by Kent and Wood (1997). All are based on the log-linearity given by the power law relationship; see (2.3) and (1.4). The simplest case is ordinary least squares [see (1.5)] for which

$$L_u = \frac{(\log u - m^{-1}\sum_{v=1}^m \log v)}{\sum_{u=1}^m (\log u - m^{-1}\sum_{v=1}^m \log v)^2}, \qquad u=1,\ldots,m.$$

Define

$$(2.6) \qquad \hat{\alpha} = \sum_{u=1}^m L_u \log \bar{Z}_u \quad \text{and} \quad \alpha_n = \sum_{u=1}^m L_u \log \mu_u$$

and

(2.7) $$F = (F_1 + F_2)/G_1,$$

where

(2.8) $$F_1 = \left\{ \sum_{u=1}^{m} L_u \mu_{0,u}^{-1} (\tau_{0,1u} + \tau_{0,2u} + \tau_{0,3u}) \right\} \\ \times \int_{t \in [0,1]^d} g^{(1)}\{X(t)\} g^{(3)}\{X(t)\} \, dt$$

and

(2.9) $$F_2 = \left\{ \sum_{u=1}^{m} L_u \mu_{0,u}^{-1} (\tau_{0,1u} + \tau_{0,2u} + \tau_{0,4u}) \right\} \\ \times \int_{t \in [0,1]^d} [g^{(2)}\{X(t)\}]^2 \, dt.$$

In the above $\mu_{0,u} = \lim_{n\to\infty} \mu_u \in (0,\infty)$ and $\tau_{0,ju} = \lim_{n\to\infty} \tau_{ju}$ ($j=1,\ldots,4$; $u=1,\ldots,m$), where $\mu_u$ is defined in (2.4) and $\tau_{1u}$, $\tau_{2u}$, $\tau_{3u}$ and $\tau_{4u}$ are defined in (5.13), (5.18) and (5.21). Note that $\mu_{0,u}$ and $\tau_{0,ju}$ may be determined explicitly using (1.1) when $d=1$ or (2.3) when $d=2$, but we omit the details because these formulas are not required in what follows.

The quantity $\sigma$ which appears in Theorem A is defined as follows. When $d=1$,

(2.10) $$\sigma = \sqrt{t^T \Phi_0 t},$$

where $t_u = C^{-1} u^{-\alpha} L_u$, $u=1,\ldots,m$, $t = (t_1,\ldots,t_m)^T$, $C > 0$ is the constant in formula (2.7) in Kent and Wood (1997) and $\Phi_0$ is the covariance matrix defined via (3.3) and (3.6) in Kent and Wood (1997) (the precise definitions of $C$ and $\Phi_0$ need not concern us here). When $d=2$, $\sigma$ is still of the form (2.10), but with the quantities $C$ and $\Phi_0$ now given by (3.13) and (3.16) in Chan and Wood (2000).

The integrator $B(\cdot)$ which appears in Theorem A is a Gaussian random measure such that, for any measurable subsets $A_1,\ldots,A_k$ of $[0,1]^d$,

(2.11) $$(B(A_1),\ldots,B(A_k))^T \sim N_k(0,\Psi),$$

where $\Psi = (\psi_{jl})$, $\psi_{jl} = \lambda(A_j \cap A_l)$ and $\lambda$ denotes Lebesgue measure on $\mathbf{R}^d$.

In Theorem B $Z_S$ is the zero-mean random Gaussian measure defined on $[-\pi,\pi]^d$ with the following properties: if $D_1, D_2 \subseteq [-\pi,\pi]^d$, then

$$\mathrm{cov}\{Z_S(D_1), Z_S(D_2)\} = S(D_1 \cap D_2),$$

where $S$ is the spectral measure of a covariance function on $\mathbf{Z}^d$ of the form $\rho(k) = \|k\|^{-\alpha} A(k/\|k\|)$, where $k \in \mathbf{Z}^d$ and $A(\cdot) > 0$ is a continuous

function on the unit sphere in $\mathbf{R}^d$. In fact, $\rho(k)$ is defined as the limit as $n \to \infty$ of the covariance $\text{cov}(Y_i, Y_{i+k})$, where $Y_i = \sum_{u=1}^m u^2 Y_{iu}/(\sum_{u=1}^m u^4)^{1/2}$ and $Y_{iu}$ is defined in (5.1). When $d = 1$, we may take $A$ to be a constant, and when $d = 2$, $A$ is determined by $M$ in (2.3).

2.2. *The theorems.* We are now ready to state our main results, Theorems A and B. Discussion of these results follows in Section 2.3.

THEOREM A. *Suppose that $\gamma$ and $g$ satisfy conditions $(\mathcal{A}1)_4^{(d)}$, $(\mathcal{A}2)$ and $(\mathcal{A}3)$, and, if $d = 2$, suppose also that the sampling regime satisfies $(\mathcal{A}4)$. Let $\hat{\alpha}$ and $\alpha_n$ be defined as in (2.6), using the index set (2.1) in (2.4). Suppose that $\mathbf{a} = \{a_j\}$ is an increment of order $p \geq 0$ such that $4 + 4p - 2\alpha > d$, where $d = 1$ or $d = 2$. Then*

$$(2.12) \qquad \hat{\alpha} - \alpha_n = n^{-\alpha/d} F + n^{-1/2} \sigma \frac{\sqrt{G_2}}{G_1} Z + o_p(n^{-\alpha/d} + n^{-1/2}),$$

*where $\sigma > 0$ is a constant and $Z \sim N(0,1)$ is independent of $G_1$, $G_2$ and $F$, where the latter are defined by (2.5) and (2.7)–(2.9).*

*When $0 < 2\alpha < d$, the $n^{-\alpha/d}$ term in (2.12) is dominant and, in this case,*

$$n^{\alpha/d}(\hat{\alpha} - \alpha_n) \xrightarrow{\mathcal{D}} F.$$

*When $d < 2\alpha < \min(4, 4 + 4p - d)$, the $n^{-1/2}$ term in (2.12) is dominant and, in this case,*

$$n^{1/2}(\hat{\alpha} - \alpha_n) \xrightarrow{\mathcal{D}} \sigma G_1^{-1} \int_{[0,1]^d} [g^{(1)}\{X(t)\}]^2 \, dB(t) \stackrel{\mathcal{D}}{=} \sigma \frac{\sqrt{G_2}}{G_1} Z,$$

*where $B(t)$ is the Gaussian random measure described in (2.11). When $2\alpha = d$, both terms contribute, and we have*

$$n^{1/2}(\hat{\alpha} - \alpha_n) \xrightarrow{\mathcal{D}} F + \sigma G_1^{-1} \int_{[0,1]^d} [g^{(1)}\{X(t)\}]^2 \, dB(t).$$

*In the above the Gaussian measure $\{B(t)\}$ is independent of the underlying Gaussian field $\{X(t) : t \in [0,1]^d\}$.*

When $4 + 4p - 2\alpha < d$, the limit distribution of $\hat{\alpha}$ may be expressed as a Wiener–Itô integral of order 2. See Dobrushin and Major (1979) and Major (1981) for further details on Wiener–Itô integrals.

THEOREM B. *Suppose that $\gamma$ and $g$ satisfy conditions $(\mathcal{A}1)_4^{(d)}$, $(\mathcal{A}2)$ and $(\mathcal{A}3)$, and (when $d = 2$) the sampling regime satisfies $(\mathcal{A}4)$. Let $\hat{\alpha}$ and $\alpha_n$ be defined as in (2.6) where, as before, $\mathbf{a} = \{a_j\}$ is an increment of order* $p = 0$, and suppose that either $d = 1$ and $3/2 < \alpha < 2$ or $d = 2$ and $1 < \alpha < 2$ (corresponding to the condition $4 + 4p - 2\alpha < d$). Then

$$n^{(2-\alpha)/d}(\hat{\alpha} - \alpha_n)$$
$$\xrightarrow{\mathcal{D}} G_1^{-1} \int_{x_1,x_2 \in [-\pi,\pi]^d} \int_{t \in [0,1]^d} I(x_1 \neq x_2) e^{it^T(x_1+x_2)}$$
$$\times [g^{(1)}\{X(t)\}]^2 \, dt \, dZ_S(x_1) \, dZ_S(x_2)$$
$$\stackrel{\mathcal{D}}{=} \int_{x_1 \in [-\pi,\pi]^d} \int_{x_2 \in [-\pi,\pi]^d} I(x_1 \neq x_2) \mathcal{F}_g(x_1 + x_2) \, dZ_S(x_1) \, dZ_S(x_2),$$

where the indicator function $I(x_1 \neq x_2)$ excludes the diagonal and

$$\mathcal{F}_g(x) = G_1^{-1} \int_{t \in [0,1]^d} e^{it^T x} [g^{(1)}\{X(t)\}]^2 \, dt$$

is the Fourier transform of the random probability measure on $[0,1]^d$ whose density with respect to Lebesgue measure is given by $[g^{(1)}\{X(t)\}]^2/G_1$. In the above, the Gaussian measure $Z_S$ is independent of the underlying Gaussian field $\{X(t) : t \in [0,1]^d\}$.

2.3. *Discussion.* A number of comments concerning Theorems A and B now follow.

1. When $g$ is affine in Theorem A, $F \equiv 0$, $\sqrt{G_2}/G_1 = 1$ and therefore the limit distribution of $n^{1/2}(\hat{\alpha} - \alpha_n)$ is $N(0, \sigma^2)$. This agrees with the central limit theorems given in Kent and Wood (1995, 1997) and Chan and Wood (2000).
2. Note that, in (2.12) $\sqrt{G_2}/G_1 \geq 1$ by the Cauchy–Schwarz inequality, with equality if $g$ is affine. Moreover, when $\alpha < d/2$ the rate of convergence of $\hat{\alpha}$ to $\alpha_n$ is of slower order than $n^{-1/2}$. Therefore, we may conclude that, from the point of view of the estimator $\hat{\alpha}$ in (2.6), the non-Gaussian case ($g$ not affine) is less favorable than the Gaussian case ($g$ affine) in the framework considered in this paper. This finding is confirmed by the numerical MSE results in Section 3, though these results also suggest that the deterioration is not too severe provided that $g$ is not drastically nonaffine over the relevant domain.
3. At the borderline between Theorems A and B [i.e., when the increment used has order 0 and $\alpha = (4-d)/2$], the limit distribution is of the type given in Theorem A, but the convergence rate of $n^{-1/2}$ is modified by a logarithmic factor in $n$; compare the Constantine–Hall (1994) result when $\alpha = 3/2$. We omit the proof.
4. The bias of $\hat{\alpha}$ in the context of Theorem A depends not only on the right-hand side of (2.12), but also on $\alpha_n - \alpha$. It follows from (2.1) and (2.3) and the definitions in (2.6) that $\alpha_n - \alpha = O(n^{-\beta/d})$. It can be seen that, in Theorem A, this term makes a negligible contribution to the bias if and only if $\beta > \alpha$, and it makes a negligible contribution to the mean squared error if and only if $\beta > \min(\alpha, d/2)$. In Theorem B $\alpha_n - \alpha$ makes a negligible contribution to the MSE if $\beta > 2 - \alpha$. See Constantine and Hall (1994), Kent and Wood (1997) and Chan and Wood (2000) for further discussion of this bias term.
5. It is interesting to note that Theorems A and B do not require any moment conditions on $g\{X(t)\}$. This is a consequence of the Gaussianity of $X(t)$ and the smoothness of $g$. See also the proof of Step 1 in Section 5.
6. It is straightforward to extend the results given for $d = 1, 2$ to $d > 2$. In short, if $p$ is the order of the increment used, then we are in the situation of Theorem A if $4 + 4p - 2\alpha > d$, and we are in the situation of Theorem B if $4 + 4p - 2\alpha < d$.
7. It may be helpful to give some intuition as to why $X(t)$ is independent of the Gaussian measure $B$ in Theorem A and the Gaussian measure $Z_S$ in Theorem B. For simplicity consider the case $d = 1$ and let $X(t)$ have covariance function $\gamma(t)$ satisfying (1.1). Then as $h \to 0$ (corresponding to a fixed-domain asymptotic regime), we have (i) $W_h(t) \equiv h^{-\alpha/2}\{X(t+h) - X(t)\}$ converges in distribution to $N(0, 2c)$, and (ii) $W_h(t)$ is asymptotically independent of $X(t)$. Essentially, the increments of $B$ in Theorem A and $Z_S$ in Theorem B are linear combinations of terms of the form $W_h(t)$, while each integrand depends only on $X(t)$. Thus, the integrator and integrand are independent in the limit. In Step 6 of Section 5 we establish this asymptotic independence rigorously.
8. Note that the "diagonal" $\{x_1 = x_2\}$ is explicitly excluded from the region of integration in Theorem B. This is in line with the definition of the Wiener–Itô integral given by many authors, including Dobrushin and Major (1979), formula (1.9), Major (1981), Theorem 8.2, Nualart (1995), formula (1.13) and Arcones (1994), formula (3.9), but note that in all of the above references the exclusion of the diagonal is not made explicit in the notation. See Taqqu (1979), page 77, for helpful discussion of this point.

**3. Numerical results.** In the simulation studies described below the covariance function of the underlying stationary Gaussian field was chosen to be of the form $\gamma(t) = \exp(-c\|t\|^\alpha)$, where $\alpha \in (0, 2)$ and $c > 0$, and $\|t\|$ is the usual Euclidean norm on $\mathbf{R}^d$. We chose $c = 1$ when $d = 1$ and $c = 10$ when $d = 2$ throughout our numerical work. The data were simulated using the circulant embedding approach; see, for example, Wood and Chan (1994) and Chan and Wood (1999) for details and further references. In all cases considered, it was possible to use the algorithm in its "exact" form.

Four types of point transformations were considered:

1. Uniform: $g(x) = \Phi(x)$;
2. Exponential: $g(x) = -\log\{1 - \Phi(x)\}$;
3. Chi-squared with one degree of freedom: $g(x) = x^2$;
4. Log-normal: $g(x; \tau) = \exp(\tau x)$.

In the above $\Phi$ denotes the standard normal distribution function. Two cases of the log-normal distribution were considered, corresponding to $\tau = 1$ and $\tau = 4$. Note that $\tau = 4$ corresponds to an extremely nonaffine transformation in the relevant domain and is included as an extreme case. Each of the above transformations preserves the fractal dimension of the underlying Gaussian random field [see Hall and Roy (1994)].

Figure 1(a) shows a realization of a Gaussian process with $\alpha = 0.1$, and Figure 1(b)–(f) shows various nonaffine transformations of this realization. Figure 2 shows similar nonaffine transformations of a smoother process, with $\alpha = 1.0$. It is clear from visual inspection of the graphs that the nonaffine transformations do have a noticeable effect in both Figures 1 and 2, and in many cases the transformed processes do clearly exhibit non-Gaussian features. The log-normal transformation with $\tau = 4$ is particularly extreme, as might be expected.

Some representative numerical results for $d = 1$ are displayed in Tables 1–3 and for $d = 2$ in Table 4.

Table 1 summarizes the results of a simulation study of the performance of three estimators, $\hat{\alpha}_{\text{OLS}}^{(0)}$, $\hat{\alpha}_{\text{OLS}}^{(1)}$ and $\hat{\alpha}_{\text{GLS}}^{(1)}$, of the fractal index $\alpha$ of a Gaussian process, and several transformations of this process. The notation used for these estimators is the same as that in Kent and Wood (1997) and Chan and Wood (2000): $\hat{\alpha}_{\text{OLS}}^{(p)}$, $p = 0$ or $1$, is the ordinary least squares estimator of $\alpha$ given by (1.5), based on the log–log relationship (1.4), and using increment (A.3) when $p = 0$ and increment (A.4) when $p = 1$; and $\hat{\alpha}_{\text{GLS}}^{(1)}$ is a generalized least squares estimator of $\alpha$, again based on the log–log relationship, using increment (A.4). For smoother processes (e.g., $\alpha = 1.9$) the GLS estimator performs slightly better, in terms of the mean squared error (MSE), than the other two estimators, due to its smaller bias. For rougher processes (e.g., $\alpha = 1.0$ and $0.1$) all three estimators have very similar MSE with slightly higher standard deviation (SD) for the estimators with $p = 1$. One noticeable fact is that all estimators generally perform worse for the transformed processes. However, the deterioration is fairly mild in most cases, except in the case of the log-normal(4) transformation, where there is a substantial increase in the MSE.

In the second study we compare the asymptotic and empirical rate of decrease in variance of $\hat{\alpha}_{\text{OLS}}^{(0)}$ and $\hat{\alpha}_{\text{OLS}}^{(1)}$ as the sample size $n$ increases. The asymptotic rates of decrease are computed using the following results in the

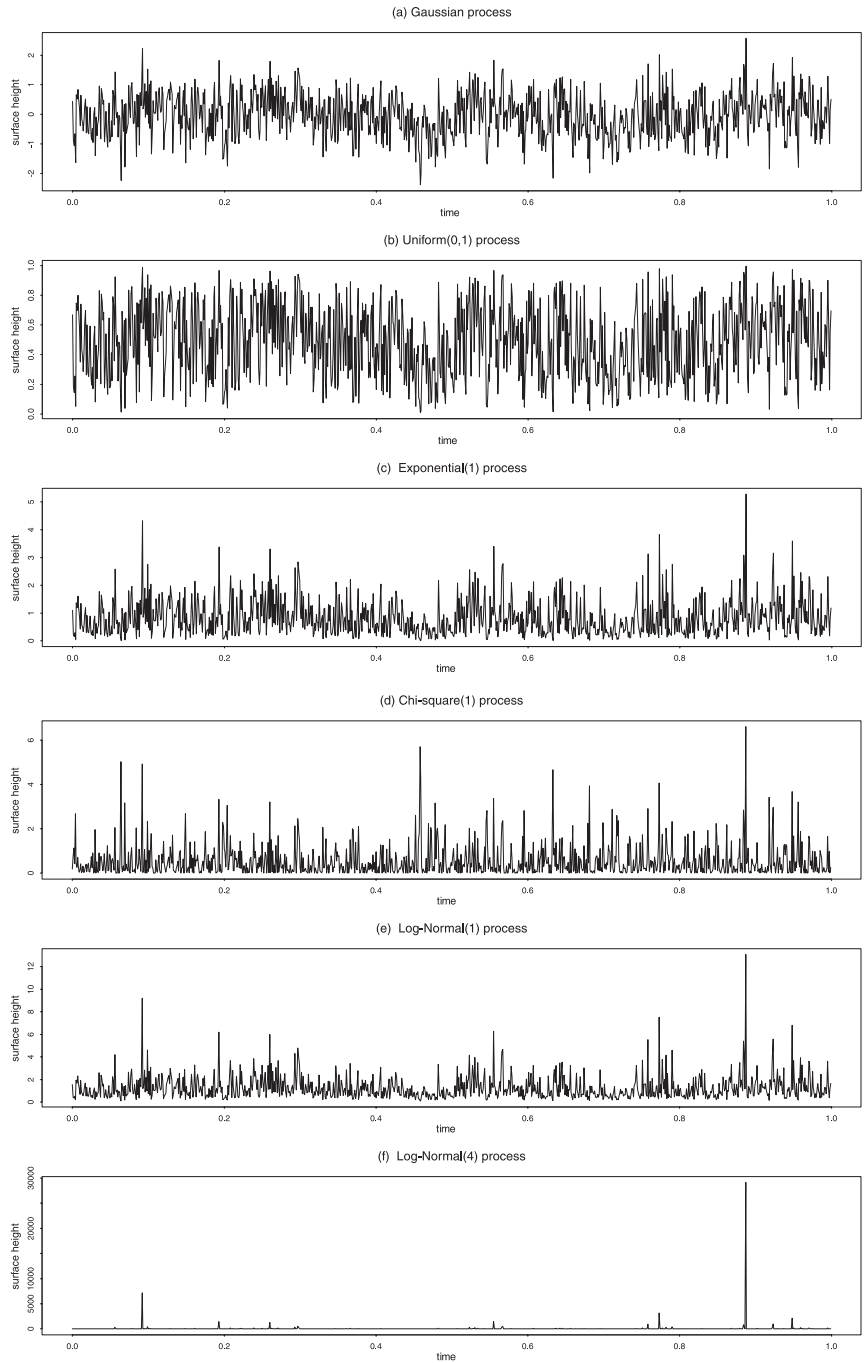

Fig. 1. *Gaussian process and its point-transformed non-Gaussian processes with $n = 1000$, $\alpha = 0.1$.*

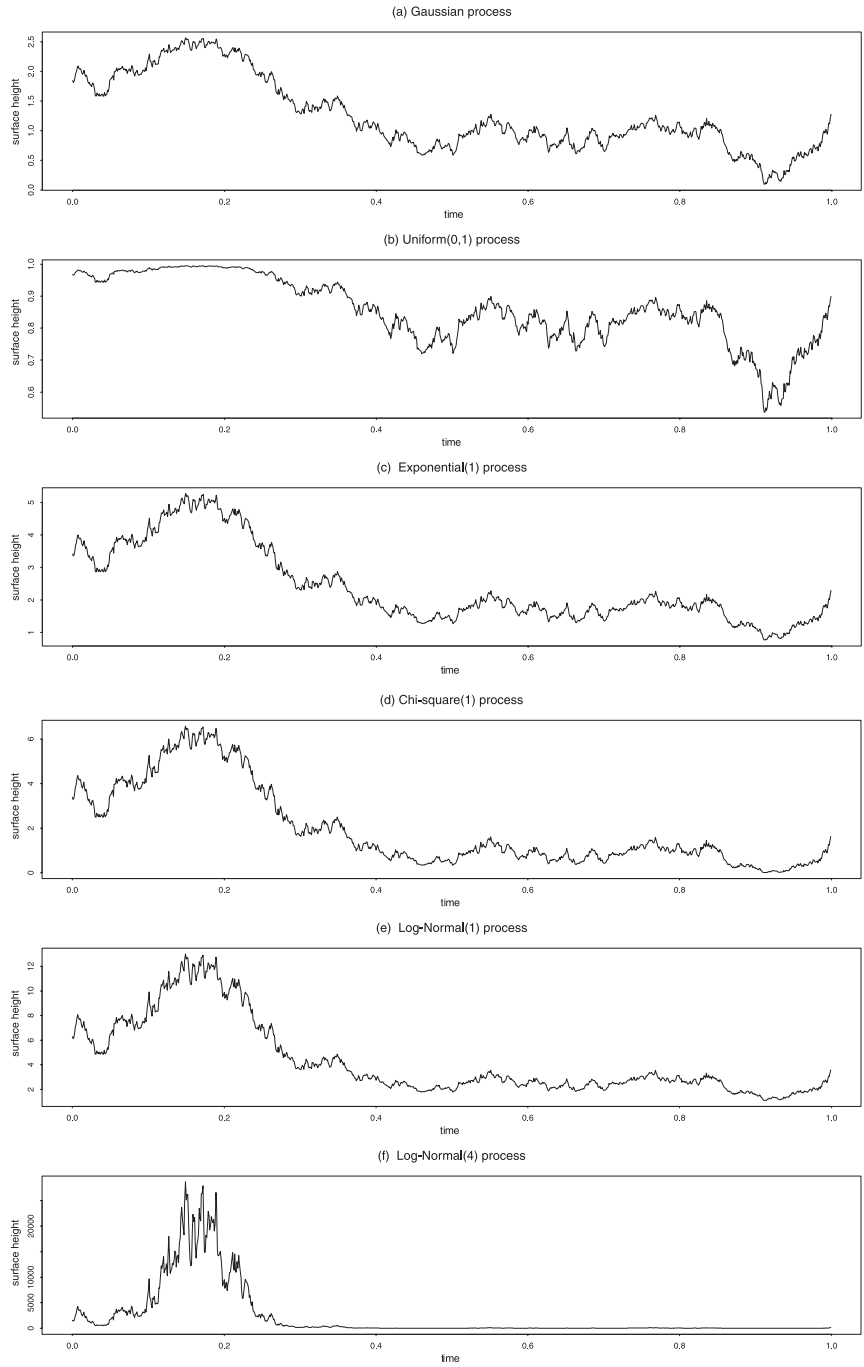

Fig. 2. *Gaussian process and its point-transformed non-Gaussian processes with $n = 1000$, $\alpha = 1.0$.*

TABLE 1
*Comparison of the three estimators and six processes with $n = 1000$, $m = 4$, based on 100 simulations in each case*

| $\alpha$ | Process | $\hat{\alpha}_{\mathrm{OLS}}^{(0)}$ | | | $\hat{\alpha}_{\mathrm{OLS}}^{(1)}$ | | | $\hat{\alpha}_{\mathrm{GLS}}^{(1)}$ | | |
|---|---|---|---|---|---|---|---|---|---|---|
| | | Bias | SD | MSE | Bias | SD | MSE | Bias | SD | MSE |
| 0.1 | Gaussian | −0.021 | 0.033 | 0.002 | −0.021 | 0.043 | 0.002 | −0.021 | 0.043 | 0.002 |
| | Uniform | −0.028 | 0.034 | 0.002 | −0.029 | 0.045 | 0.003 | −0.029 | 0.045 | 0.003 |
| | Exp(1) | −0.027 | 0.039 | 0.002 | −0.025 | 0.051 | 0.003 | −0.025 | 0.051 | 0.003 |
| | $\chi_1^2$ | −0.041 | 0.041 | 0.003 | −0.036 | 0.054 | 0.004 | −0.036 | 0.054 | 0.004 |
| | Log-N(1) | −0.032 | 0.047 | 0.003 | −0.029 | 0.064 | 0.005 | −0.029 | 0.063 | 0.005 |
| | Log-N(4) | −0.079 | 0.090 | 0.014 | −0.081 | 0.152 | 0.030 | −0.075 | 0.133 | 0.024 |
| 1.0 | Gaussian | −0.002 | 0.041 | 0.002 | 0.002 | 0.059 | 0.003 | 0.001 | 0.057 | 0.003 |
| | Uniform | −0.000 | 0.052 | 0.003 | 0.006 | 0.077 | 0.006 | 0.005 | 0.074 | 0.006 |
| | Exp(1) | −0.005 | 0.055 | 0.003 | −0.005 | 0.076 | 0.006 | −0.004 | 0.074 | 0.006 |
| | $\chi_1^2$ | −0.008 | 0.059 | 0.004 | 0.000 | 0.079 | 0.006 | 0.001 | 0.074 | 0.005 |
| | Log-N(1) | −0.008 | 0.057 | 0.003 | −0.008 | 0.079 | 0.006 | −0.006 | 0.077 | 0.006 |
| | Log-N(4) | −0.059 | 0.134 | 0.021 | −0.054 | 0.186 | 0.038 | −0.051 | 0.175 | 0.033 |
| 1.9 | Gaussian | −0.030 | 0.055 | 0.004 | −0.002 | 0.056 | 0.003 | −0.002 | 0.054 | 0.003 |
| | Uniform | −0.025 | 0.060 | 0.004 | −0.001 | 0.068 | 0.005 | −0.000 | 0.064 | 0.004 |
| | Exp(1) | −0.033 | 0.056 | 0.004 | −0.003 | 0.068 | 0.005 | −0.002 | 0.064 | 0.004 |
| | $\chi_1^2$ | −0.041 | 0.055 | 0.005 | −0.010 | 0.071 | 0.005 | −0.009 | 0.066 | 0.004 |
| | Log-N(1) | −0.038 | 0.054 | 0.004 | −0.004 | 0.072 | 0.005 | −0.003 | 0.067 | 0.005 |
| | Log-N(4) | −0.068 | 0.063 | 0.009 | −0.013 | 0.123 | 0.015 | −0.012 | 0.119 | 0.014 |

TABLE 2
*Comparison of the two OLS estimators and three processes with $n = 1000, 2000, 4000, 10,000$, $m = 10$, based on $500$ simulations in each case; the empirical (asymptotic) variance ratios are given when $n \geq 2000$*

| Process | $\alpha$ | Order 0 var | 2000 | | 4000 | | 10,000 | | Order 1 var | 2000 | | 4000 | | 10,000 | |
|---|---|---|---|---|---|---|---|---|---|---|---|---|---|---|---|
| Gaussian | 0.1 | $0.0^340$ | 0.48 | (0.50) | 0.24 | (0.25) | 0.10 | (0.10) | $0.0^360$ | 0.48 | (0.50) | 0.25 | (0.25) | 0.10 | (0.10) |
| | 0.4 | $0.0^215$ | 0.50 | (0.50) | 0.23 | (0.25) | 0.10 | (0.10) | $0.0^221$ | 0.52 | (0.50) | 0.22 | (0.25) | 0.10 | (0.10) |
| | 0.7 | $0.0^222$ | 0.47 | (0.50) | 0.23 | (0.25) | 0.09 | (0.10) | $0.0^232$ | 0.45 | (0.50) | 0.23 | (0.25) | 0.10 | (0.10) |
| | 1.0 | $0.0^225$ | 0.47 | (0.50) | 0.22 | (0.25) | 0.09 | (0.10) | $0.0^237$ | 0.49 | (0.50) | 0.25 | (0.25) | 0.10 | (0.10) |
| | 1.3 | $0.0^227$ | 0.53 | (0.50) | 0.29 | (0.25) | 0.10 | (0.10) | $0.0^236$ | 0.55 | (0.50) | 0.29 | (0.25) | 0.11 | (0.10) |
| | 1.6 | $0.0^237$ | 0.69 | (0.57) | 0.33 | (0.33) | 0.19 | (0.16) | $0.0^239$ | 0.52 | (0.50) | 0.27 | (0.25) | 0.11 | (0.10) |
| | 1.9 | $0.0^241$ | 0.73 | (0.87) | 0.62 | (0.76) | 0.40 | (0.63) | $0.0^246$ | 0.44 | (0.50) | 0.22 | (0.25) | 0.09 | (0.10) |
| Exp(1) | 0.1 | $0.0^349$ | 0.56 | (0.87) | 0.28 | (0.76) | 0.14 | (0.63) | $0.0^373$ | 0.56 | (0.87) | 0.28 | (0.76) | 0.14 | (0.63) |
| | 0.4 | $0.0^223$ | 0.57 | (0.57) | 0.30 | (0.33) | 0.12 | (0.16) | $0.0^234$ | 0.53 | (0.57) | 0.28 | (0.33) | 0.11 | (0.16) |
| | 0.7 | $0.0^233$ | 0.50 | (0.50) | 0.27 | (0.25) | 0.11 | (0.10) | $0.0^247$ | 0.50 | (0.50) | 0.28 | (0.25) | 0.12 | (0.10) |
| | 1.0 | $0.0^234$ | 0.58 | (0.50) | 0.27 | (0.25) | 0.11 | (0.10) | $0.0^258$ | 0.55 | (0.50) | 0.26 | (0.25) | 0.10 | (0.10) |
| | 1.3 | $0.0^238$ | 0.53 | (0.50) | 0.34 | (0.25) | 0.11 | (0.10) | $0.0^253$ | 0.58 | (0.50) | 0.30 | (0.25) | 0.11 | (0.10) |
| | 1.6 | $0.0^245$ | 0.71 | (0.57) | 0.32 | (0.33) | 0.20 | (0.16) | $0.0^261$ | 0.50 | (0.50) | 0.28 | (0.25) | 0.10 | (0.10) |
| | 1.9 | $0.0^242$ | 0.74 | (0.87) | 0.63 | (0.76) | 0.41 | (0.63) | $0.0^258$ | 0.51 | (0.50) | 0.25 | (0.25) | 0.11 | (0.10) |
| Log-normal(4) | 0.1 | $0.0^212$ | 0.75 | (0.87) | 0.91 | (0.76) | 0.57 | (0.63) | $0.0^233$ | 0.69 | (0.87) | 0.80 | (0.76) | 0.53 | (0.63) |
| | 0.4 | $0.0^111$ | 1.0 | (0.57) | 0.76 | (0.33) | 0.91 | (0.16) | $0.0^121$ | 0.92 | (0.57) | 0.62 | (0.33) | 0.74 | (0.16) |
| | 0.7 | $0.0^118$ | 0.73 | (0.50) | 0.54 | (0.25) | 0.36 | (0.10) | $0.0^136$ | 0.77 | (0.50) | 0.49 | (0.25) | 0.31 | (0.10) |
| | 1.0 | $0.0^117$ | 0.63 | (0.50) | 0.35 | (0.25) | 0.18 | (0.10) | $0.0^137$ | 0.56 | (0.50) | 0.32 | (0.25) | 0.14 | (0.10) |
| | 1.3 | $0.0^111$ | 0.52 | (0.50) | 0.38 | (0.25) | 0.15 | (0.10) | $0.0^129$ | 0.50 | (0.50) | 0.27 | (0.25) | 0.11 | (0.10) |
| | 1.6 | $0.0^283$ | 0.70 | (0.57) | 0.36 | (0.33) | 0.22 | (0.16) | $0.0^124$ | 0.49 | (0.50) | 0.25 | (0.25) | 0.09 | (0.10) |
| | 1.9 | $0.0^246$ | 0.77 | (0.87) | 0.69 | (0.76) | 0.50 | (0.63) | $0.0^118$ | 0.48 | (0.50) | 0.26 | (0.25) | 0.11 | (0.10) |

Gaussian case [see Constantine and Hall (1994) and Chan and Wood (2000)]:

$$\mathrm{var}(\hat{\alpha}_{\mathrm{OLS}}^{(0)}) \sim \begin{cases} C_1 n^{-1}, & 0 < \alpha < 3/2, \\ C_2 n^{-1} \log n, & \alpha = 3/2, \\ C_3 n^{2\alpha-4}, & 3/2 < \alpha < 2, \end{cases}$$

and

$$\mathrm{var}(\hat{\alpha}_{\mathrm{OLS}}^{(1)}) = C_4 n^{-1}, \qquad 0 < \alpha < 2,$$

where $C_1, \ldots, C_4$ depend on $\alpha$ and $m$ but not on $n$. Hence, the asymptotic variance ratios for sample sizes $n_1 < n_2$ and $\alpha \neq 3/2$ will be equal to either $(n_2/n_1)^{-1}$ or $(n_2/n_1)^{2\alpha-4}$. For a non-Gaussian process in the one-dimensional case ($d=1$), the variance formulas for both $\hat{\alpha}_{\mathrm{OLS}}^{(0)}$ and $\hat{\alpha}_{\mathrm{OLS}}^{(1)}$ are the same as in the Gaussian case except when $0 < \alpha < 1/2$. In that case the asymptotic variances are asymptotic to $C_5 n^{-2\alpha}$ for both estimators. Table 2 reports the estimated variance when $n = 1000$ and the estimated rate of decrease in variance:

$$\widehat{\mathrm{var}}(\hat{\alpha}_n)/\widehat{\mathrm{var}}(\hat{\alpha}_{1000}), \qquad \text{when } n = 2000, 4000, 10{,}000,$$

for Gaussian, exponential(1) (as an example of mildly nonaffine transformation) and log-normal(4) (as an example of extremely nonaffine transformation) processes. For the Gaussian case there is good agreement between the theoretical and the numerical results. For the non-Gaussian cases, there is good agreement for medium to large $\alpha$ (smoother process). For the exponential(1) process, the numerical variance ratios are closer to those for the Gaussian process for small $\alpha$ (rough process). The rate of decrease in variance does not depend noticeably on the choice of $m$.

In Table 3 we study the effect of the number of points $m$ used in the regression for $\hat{\alpha}_{\mathrm{OLS}}^{(0)}$, $\hat{\alpha}_{\mathrm{OLS}}^{(1)}$ and $\hat{\alpha}_{\mathrm{GLS}}^{(1)}$. Table 3 suggests that the number of points used in the regression does not affect the MSE significantly for all three estimators. For these simulated data a choice of $m = 4$ would be suitable.

Next we take a closer look at the distribution of $\hat{\alpha}_{\mathrm{GLS}}^{(1)}$ for both Gaussian and non-Gaussian processes. Figure 3 shows normal quantile–quantile plots based on 100 estimates for Gaussian, chi-square(1) and log-normal(4) processes. The added straight lines go through the first and third quartiles of these estimates and the corresponding value of the standard normal distribution. With the exception of the top right-hand corner [log-normal(4) with $\alpha = 0.1$], all cases suggest that these empirical distributions are reasonably close to normal. This graphical finding is supported by the Kolmogorov–Smirnov goodness-of-fit test.

Our tentative conclusions for $d=1$ are as follows:

1. All three estimators perform fairly well under modest departures from normality (of the type introduced by the nonaffine transformation $g$).
2. Our numerical results suggest there is no advantage in using the GLS estimator in the non-Gaussian case (which is not surprising, as it was designed for the Gaussian case).
3. The number of points used in the regression is not critical and can be taken as small as 4 for the simulation data considered.

For $d=2$, the performance of the OLS estimator, based on the "square" increment defined in (A.8), is studied for $n_0 = (50, 50), (100, 100)$ and $(500, 500)$, and $m = 2(1)10$. We denote this estimator by $\hat{\alpha}_{n_0}$. In particular, we compare the asymptotic and empirical rates of decrease in variance of $\hat{\alpha}_{n_0}$ as the sampling region increases. Theoretical results from Chan and Wood (2000) imply that for the Gaussian case,

$$\text{var}(\hat{\alpha}_{(n_1, n_1)})/\text{var}(\hat{\alpha}_{(n_2, n_2)}) = (n_2/n_1)^2, \qquad 0 < \alpha < 2.$$

Theorem A implies that for non-Gaussian fields,

$$\text{var}(\hat{\alpha}_{(n_1, n_1)})/\text{var}(\hat{\alpha}_{(n_2, n_2)}) = \begin{cases} (n_2/n_1)^2, & 1 < \alpha < 2, \\ (n_2/n_1)^{2\alpha}, & 0 < \alpha < 1. \end{cases}$$

TABLE 3
*Comparison of the three estimators for the $\chi_1^2$ process in terms of MSE with $n = 2000$, $m = 2, 4, 6, 8, 10$, based on 100 simulations in each case*

| | | $m$ | | | | |
|---|---|---|---|---|---|---|
| $\alpha$ | $\hat{\alpha}$ | 2 | 4 | 6 | 8 | 10 |
| 0.1 | $\hat{\alpha}_{\text{OLS}}^{(0)}$ | 0.0056 | 0.0031 | 0.0032 | 0.0031 | 0.0031 |
| | $\hat{\alpha}_{\text{OLS}}^{(1)}$ | 0.0103 | 0.0035 | 0.0033 | 0.0030 | 0.0031 |
| | $\hat{\alpha}_{\text{GLS}}^{(1)}$ | 0.0103 | 0.0035 | 0.0033 | 0.0030 | 0.0031 |
| 0.3 | $\hat{\alpha}_{\text{OLS}}^{(0)}$ | 0.0054 | 0.0037 | 0.0041 | 0.0045 | 0.0048 |
| | $\hat{\alpha}_{\text{OLS}}^{(1)}$ | 0.0111 | 0.0041 | 0.0040 | 0.0043 | 0.0045 |
| | $\hat{\alpha}_{\text{GLS}}^{(1)}$ | 0.0111 | 0.0041 | 0.0038 | 0.0039 | 0.0040 |
| 1.0 | $\hat{\alpha}_{\text{OLS}}^{(0)}$ | 0.0017 | 0.0018 | 0.0019 | 0.0020 | 0.0022 |
| | $\hat{\alpha}_{\text{OLS}}^{(1)}$ | 0.0053 | 0.0038 | 0.0039 | 0.0035 | 0.0033 |
| | $\hat{\alpha}_{\text{GLS}}^{(1)}$ | 0.0053 | 0.0035 | 0.0028 | 0.0023 | 0.0021 |
| 1.7 | $\hat{\alpha}_{\text{OLS}}^{(0)}$ | 0.0021 | 0.0025 | 0.0029 | 0.0032 | 0.0035 |
| | $\hat{\alpha}_{\text{OLS}}^{(1)}$ | 0.0054 | 0.0033 | 0.0029 | 0.0032 | 0.0036 |
| | $\hat{\alpha}_{\text{GLS}}^{(1)}$ | 0.0054 | 0.0030 | 0.0027 | 0.0024 | 0.0022 |
| 1.9 | $\hat{\alpha}_{\text{OLS}}^{(0)}$ | 0.0028 | 0.0031 | 0.0033 | 0.0035 | 0.0037 |
| | $\hat{\alpha}_{\text{OLS}}^{(1)}$ | 0.0048 | 0.0032 | 0.0032 | 0.0034 | 0.0037 |
| | $\hat{\alpha}_{\text{GLS}}^{(1)}$ | 0.0048 | 0.0027 | 0.0024 | 0.0022 | 0.0023 |

TABLE 4
*Comparison between empirical and asymptotic variance ratios among Gaussian and non-Gaussian random fields with $n_0 = (50, 50)$, $(100, 100)$, $(500, 500)$, $m = 4$, based on $100$ simulations in each case*

| $\alpha$ | var | 100/50 | | 500/50 | | 500/100 | | var | 100/50 | | 500/50 | | 500/100 | |
|---|---|---|---|---|---|---|---|---|---|---|---|---|---|---|
| | | Gaussian field | | | | | | | Uniform field | | | | | |
| 0.1 | $0.0^2 14$ | 0.15 | (0.25) | 0.01 | (0.01) | 0.04 | (0.04) | $0.0^2 13$ | 0.15 | (0.87) | 0.01 | (0.63) | 0.04 | (0.72) |
| 0.4 | $0.0^2 15$ | 0.23 | (0.25) | 0.01 | (0.01) | 0.04 | (0.04) | $0.0^2 15$ | 0.21 | (0.57) | 0.01 | (0.16) | 0.05 | (0.28) |
| 0.7 | $0.0^2 30$ | 0.24 | (0.25) | 0.01 | (0.01) | 0.03 | (0.04) | $0.0^2 29$ | 0.22 | (0.38) | 0.01 | (0.04) | 0.05 | (0.11) |
| 1.0 | $0.0^2 41$ | 0.20 | (0.25) | 0.01 | (0.01) | 0.04 | (0.04) | $0.0^2 52$ | 0.24 | (0.25) | 0.01 | (0.01) | 0.04 | (0.04) |
| 1.3 | $0.0^2 49$ | 0.21 | (0.25) | 0.01 | (0.01) | 0.03 | (0.04) | $0.0^2 69$ | 0.28 | (0.25) | 0.01 | (0.01) | 0.03 | (0.04) |
| 1.6 | $0.0^2 55$ | 0.22 | (0.25) | 0.01 | (0.01) | 0.05 | (0.04) | $0.0^1 11$ | 0.34 | (0.25) | 0.01 | (0.01) | 0.03 | (0.04) |
| 1.9 | $0.0^1 12$ | 0.26 | (0.25) | 0.01 | (0.01) | 0.02 | (0.04) | $0.0^1 38$ | 0.48 | (0.25) | 0.01 | (0.01) | 0.03 | (0.04) |
| | | $\chi_1^2$ field | | | | | | | Log-normal(4) field | | | | | |
| 0.1 | $0.0^2 15$ | 0.23 | (0.87) | 0.01 | (0.63) | 0.03 | (0.72) | 2.98 | 0.23 | (0.87) | 0.01 | (0.63) | 0.03 | (0.72) |
| 0.4 | $0.0^2 20$ | 0.21 | (0.57) | 0.01 | (0.16) | 0.05 | (0.28) | 1.95 | 0.32 | (0.57) | 0.01 | (0.16) | 0.03 | (0.28) |
| 0.7 | $0.0^2 49$ | 0.29 | (0.38) | 0.01 | (0.04) | 0.05 | (0.11) | 3.59 | 0.06 | (0.38) | 0.01 | (0.04) | 0.01 | (0.11) |
| 1.0 | $0.0^2 97$ | 0.23 | (0.25) | 0.01 | (0.01) | 0.04 | (0.04) | 6.60 | 0.23 | (0.25) | 0.01 | (0.01) | 0.04 | (0.04) |
| 1.3 | $0.0^1 13$ | 0.25 | (0.25) | 0.01 | (0.01) | 0.04 | (0.04) | 1.36 | 1.20 | (0.25) | 0.03 | (0.01) | 0.03 | (0.04) |
| 1.6 | $0.0^1 16$ | 0.21 | (0.25) | 0.01 | (0.01) | 0.05 | (0.04) | 1.00 | 0.54 | (0.25) | 0.04 | (0.01) | 0.02 | (0.04) |
| 1.9 | $0.0^1 26$ | 0.33 | (0.25) | 0.01 | (0.01) | 0.02 | (0.04) | 2.87 | 0.11 | (0.25) | 0.04 | (0.04) | 0.41 | (0.04) |

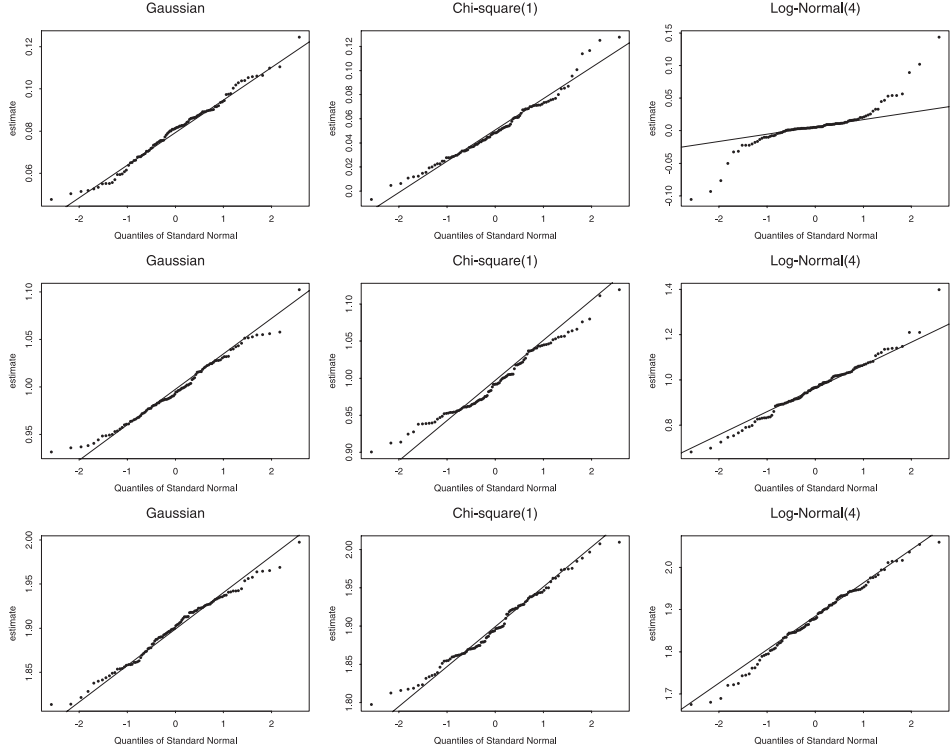

FIG. 3. *Normal quantile–quantile plots for $\hat{\alpha}^{(1)}_{\mathrm{GLS}}$ with $n = 2000$, $\alpha = 0.1, 1.0, 1.9$, $m = 10$.*

Table 4 reports the estimated variance when $n_0 = (50, 50)$ and the estimated variance ratios when $n_1 > n_2$ for $n_1, n_2 = 50, 100, 500$, for the following types of fields: Gaussian, uniform, $\chi_1^2$ and log-normal(4). For the Gaussian case there is good agreement between the theoretical and numerical results. For non-Gaussian fields there is close agreement for medium to large $\alpha$. For small $\alpha$ the empirical ratios are closer to those in the Gaussian case. As in the one-dimensional case, the number of points used in the regression is not critical. In summary, the OLS estimator with $p = 1$ performs reasonably well under mild departures from Gaussianity.

**4. Some lemmas.** We now present some results which are required in Section 5. Lemmas 4.1 and 4.2 are used to prove Lemma 4.3 which (along with the elementary Lemma 4.4) plays a key role in bounding remainder terms which arise in Steps 2–5 in Section 5. Lemma 4.5 is used in Step 6 of Section 5.

LEMMA 4.1. *Let $(A_i^T, B_i)^T$, $i \in \mathbf{Z}^d$, where $A_i = (A_{i1}, \ldots, A_{iK})^T$, be a stationary Gaussian vector field and assume that $(A_i^T, B_i)^T$ has a standard* multivariate Gaussian distribution (i.e., with mean the zero vector and identity covariance matrix). Define

$$\sigma_{AA}(i-j) = \max_{k,l=1,\ldots,K} |\text{cov}(A_{ik}, A_{jl})|, \qquad \sigma_{BB}(i-j) = |\text{cov}(B_i, B_j)|,$$

$$\sigma_{AB}(i-j) = \max_{k=1,\ldots,K} |\text{cov}(A_{ik}, B_j)|$$
$$= \max_{k=1,\ldots,K} |\text{cov}(B_i, A_{jk})| = \sigma_{BA}(i-j).$$

Let $a$ denote a nonnegative multi-index $a = (a[1], \ldots, a[K]) \in \mathbf{Z}^K$, and write $H_a(A_i) = \prod_{k=1}^{K} H_{a[k]}(A_{ik})$, where $H_m$, $m \geq 0$, is the Hermite polynomial of degree $m$. Then

$$|E[H_a(A_i) H_m(B_i) H_a(A_j) H_m(B_j)]|$$
$$\leq (a_0! m!)^2 \sum_{r=\max(a_0-m,0)}^{a_0} \frac{\sigma_{AA}(i-j)^r \sigma_{AB}(i-j)^{2(a_0-r)}}{r!(a_0-r)!(a_0-r)!(m-a_0+r)!},$$

where $a_0 = |a| = \sum_{k=1}^{K} a[k]$.

PROOF. From the diagram formula for moments (see Remark 4.1) we have, under the assumptions of the lemma,

$$E[H_a(A_i) H_m(B_i) H_a(A_j) H_m(B_j)]$$
(4.1)
$$= \frac{(a[1]! \cdots a[K]! m!)^2}{q!} \sum \sigma_{l_1,t_1}(i-j) \cdots \sigma_{l_q,t_q}(i-j),$$

where $q = m + a_0$, the $l$'s are associated with the components of $(A_i^T, B_i)^T$, the $t$'s are associated with the components of $(A_j^T, B_j)^T$ and the summation is over all indices $l_1, t_1, \ldots, l_q, t_q \in \{1, \ldots, K+1\}$ such that there are precisely $a[k]$ $l$-indices and $t$-indices equal to $k = 1, \ldots, K$, and $m$ $l$-indices and $t$-indices equal to $K+1$. In the above, $\sigma_{r,s}(i-j) = \text{cov}(A_{ir}, A_{js})$ for $1 \leq r$, $s \leq K$, $\sigma_{r,K+1}(i-j) = \text{cov}(A_{ir}, B_j)$ for $1 \leq r \leq K$, $\sigma_{K+1,s}(i-j) = \text{cov}(B_i, A_{js})$ and $\sigma_{K+1,K+1}(i-j) = \text{cov}(B_i, B_j)$.

Consider a typical product $\sigma_{l_1,k_1}(i-j) \cdots \sigma_{l_q,k_q}(i-j)$. If this consists of $r$ pairings of components of $A_i$ with components of $A_j$, then there must be $a_0 - r$ pairings of components of $A_i$ with $B_j$, $a_0 - r$ pairings of components of $A_j$ with $B_i$ and $m - a_0 + r$ pairings of $B_i$ with $B_j$, where necessarily $\max(a_0 - m, 0) \leq r \leq a_0$. Therefore,

$$|\sigma_{l_1,k_1}(i-j) \cdots \sigma_{l_q,k_q}(i-j)|$$
$$\leq \sigma_{AA}(i-j)^r \sigma_{AB}(i-j)^{a_0-r} \sigma_{BA}(i-j)^{a_0-r} \sigma_{BB}(i-j)^{m-a_0+r}$$
$$\leq \sigma_{AA}(i-j)^r \sigma_{AB}(i-j)^{2(a_0-r)},$$

since $0 \leq \sigma_{BB}(i-j) \leq 1$ and $\sigma_{AB}(i-j) = \sigma_{BA}(i-j)$. But the number of terms in the sum in (4.1) with precisely $r$ pairings of components of $A_i$ with components of $A_j$ is bounded above by

$$\left(\frac{a_0!}{\prod_{k=1}^K a_k!}\right)^2 \left(\frac{(m+a_0)!}{r!(a_0-r)!(a_0-r)!(m-a_0+r)!}\right).$$

Therefore, since $q = m + a_0$, we have the bound

$$(4.2) \quad \left|\sum \sigma_{l_1,k_1}(i-j) \cdots \sigma_{l_q,k_q}(i-j)\right|$$
$$\leq \left(\frac{a_0!}{\prod_{k=1}^K a_k!}\right)^2 \sum_{r=\max(a_0-m,0)}^{a_0} \frac{\sigma_{AA}(i-j)^r \sigma_{AB}(i-j)^{2(a_0-r)} q!}{r!(a_0-r)!(a_0-r)!(m-a_0+r)!},$$

and the lemma follows after application of the bound (4.2) to the right-hand side of (4.1). □

REMARK 4.1. For details of the diagram formula for moments of products of Hermite polynomials in Gaussian variables, see, for example, Taqqu (1977), Major (1981) and Arcones (1994). Note that the expectation on the left-hand side of (4.1) reduces to the expression on the right-hand side of (4.1) because the components of $(A_i^T, B_i)^T$ are independent Gaussian variables for each $i$. Note also that the factor $2^{-q}$ which appears in Major's (1981) version of the formula does not appear in (4.1) because we have employed the convention that, for each $k$, $l_k$ and $t_k$ in $\sigma_{l_k,t_k}(i-j)$ are such that $l_k$ is always associated with an $i$-index and $t_k$ is always associated with a $j$-index.

LEMMA 4.2. *Let* $f: \mathbf{R} \to \mathbf{R}$ *denote a function with compact support whose $q$th derivative $f^{(q)}$ is continuous on $\mathbf{R}$. Write $H_m(x)$ for the $m$th Hermite polynomial and $\phi(x)$ for the standard normal density, and let $c_m = \int_{-\infty}^{\infty} f(x) H_m(x) \phi(x) \, dx$ denote the $m$th coefficient in the expansion of $f$ in Hermite polynomials. Then $\sum_{m=0}^{\infty} c_{m+q}^2 / m! < \infty$.*

PROOF. By assumption, $f^{(q)}$ is continuous on $\mathbf{R}$ and has compact support. Therefore, $f^{(q)}$ has an expansion in Hermite polynomials of the form $\sum_{m=0}^{\infty} c_m^{(q)} \times H_m(x)/m!$, which is $L^2$-convergent in the sense that $\sum_{m=0}^{\infty} (c_m^{(q)})^2 / m! < \infty$. But repeated integration by parts using the identity $\int_{-\infty}^y \phi(x) H_m(x) \, dx = -\phi(y) \times H_{m-1}(y)$ for $m = 1, 2, \ldots$ shows that $c_m^{(q)} = c_{m+q}$, which proves the lemma. □

Lemmas 4.1 and 4.2 are used to prove the following result.

22 G. CHAN AND A. T. A. WOOD

LEMMA 4.3. *Let* $(A_i^{(n)T}, B_i^{(n)})^T$, $i \in \mathbf{Z}^d$, $n = 1, 2, \ldots$, *be a sequence of stationary Gaussian vector fields, where* $A_i^{(n)} = (A_{i1}^{(n)}, \ldots, A_{iK}^{(n)})^T$ *is a zero-mean Gaussian vector whose dimension does not depend on* $n$, *and* $B_i^{(n)} \sim N(0,1)$. *Suppose that* (i) *for each* $n$ $B_i^{(n)}$ *is independent of* $A_i^{(n)}$ *(but not necessarily independent of* $A_j^{(n)}$ *when* $i \neq j$) *and* (ii) *the smallest eigenvalue of* $\mathrm{cov}(A_i^{(n)})$ *is bounded away from 0 as* $n \to \infty$. *Let* $\pi(A_i^{(n)})$ *be a polynomial of degree* $q$ *in the components of* $A_i^{(n)}$ *such that* $E[\pi(A_i^{(n)})] = 0$ *for all* $n$. *Suppose that*

$$\sigma_{AA}^{(n)}(i-j) = \sup_{k,l=1,\ldots,K} |\mathrm{cov}(A_{ik}^{(n)}, A_{jl}^{(n)})| \leq C\{1 + |i-j|\}^{\alpha - 2p - 2}$$

*and*

$$\sigma_{AB}^{(n)}(i-j) = \sup_{k=1,\ldots,K} |\mathrm{cov}(A_{ik}^{(n)}, B_j^{(n)})| \leq C n^{-\alpha/(2d)}\{1 + |i-j|\}^{\alpha - p - 1},$$

*where* $p$ *is a nonnegative integer and* $C > 0$ *and* $\alpha \in (0,2)$ *are constants independent of* $n$. *Let* $h:\mathbf{R} \to \mathbf{R}$ *be expressible as a sum of the form* $h = h_1 + h_2$, *where* $h_1$ *is a polynomial and* $h_2$ *is a function of compact support whose* $q$th *derivative is square integrable over* $\mathbf{R}$. *Then for* $d = 1, 2$,

$$\mathrm{var}\left(n^{-1} \sum_{i \in \mathcal{I}_n} \pi(A_i^{(n)}) h(B_i^{(n)})\right) = \begin{cases} O(n^{-1}), & \text{if } p = 0 \text{ and } (2 - \alpha) > d, \\ O(n^{-1} \log n), & \text{if } p = 0 \text{ and } (2 - \alpha) = d, \\ O(n^{(\alpha - 2)/d}), & \text{if } p = 0 \text{ and } (2 - \alpha) < d, \\ O(n^{-1}), & \text{if } p \geq 1. \end{cases}$$

PROOF. Assumption (ii), combined with the assumption that the elements of $\mathrm{cov}(A_i^{(n)})$ are bounded above by $C$, implies that we may without loss of generality assume that $A_i^{(n)}$ is a standard multivariate normal vector for each $n$. Then, using multi-index notation again, we may write the polynomial $\pi(A_i^{(n)})$ as a sum of the form $\sum_a c_a H_a(A_i^{(n)})$, where $a$ ranges over a finite set of multi-indices in $\mathbf{Z}^K$ and $c_a \in \mathbf{R}$. Therefore, the result will follow for a general polynomial $\pi$ of degree $q$ if we can prove that it holds for each product $H_a(A_i^{(n)}) = \prod_k H_{a[k]}(A_{ik}^{(n)})$ of degree at most $q$.

Let the Hermite polynomial expansion of $h$ (which is convergent in the $L^2$ sense) be given by $h(x) = \sum_{m=0}^{\infty} b_m H_m(x)/m!$. Using Lemma 4.1, we obtain

$$\mathrm{var}\left(n^{-1} \sum_{i \in \mathcal{I}_n} H_a(A_i^{(n)}) h(B_i^{(n)})\right)$$

$$= \mathrm{var}\left(n^{-1} \sum_{i \in \mathcal{I}_n} H_a(A_i^{(n)}) \sum_{m=0}^{\infty} \frac{b_m}{m!} H_m(B_i^{(n)})\right)$$

$$\leq n^{-2} \sum_{m=0}^{\infty} \left(\frac{b_m}{m!}\right)^2 \left[\sum_{i,j\in\mathcal{I}_n} \left|E[H_a(A_i^{(n)})H_a(A_j^{(n)})H_m(B_i^{(n)})H_m(B_j^{(n)})]\right|\right]$$

$$\leq n^{-2} \sum_{m=0}^{\infty} \left(\frac{b_m}{m!}\right)^2$$

$$\times \left[\sum_{i,j\in\mathcal{I}_n} (a_0!m!)^2 \right.$$

$$\left. \times \sum_{r=\max(a_0-m,0)}^{a_0} \frac{\sigma_{AA}^{(n)}(i-j)^r \sigma_{AB}^{(n)}(i-j)^{2(a_0-r)}}{r!(a_0-r)!(a_0-r)!(m-a_0+r)!}\right]$$

$$= (a_0!)^2 \sum_{m=0}^{\infty} b_m^2 \left[\sum_{r=\max(a_0-m,0)}^{a_0} \{r!(a_0-r)!(a_0-r)!(m-a_0+r)!\}^{-1}\right.$$

$$\left. \times n^{-2} \sum_{i,j\in\mathcal{I}_n} \sigma_{AA}^{(n)}(i-j)^r \sigma_{AB}^{(n)}(i-j)^{2(a_0-r)}\right].$$

Using the elementary result [see Chan and Wood (2000), page 364] that

$$n^{-2} \sum_{i,j\in\mathcal{I}_n} (1+|i-j|)^{-\rho} = \begin{cases} O(n^{-1}), & \text{if } \rho > d, \\ O(n^{-1}L(n)), & \text{if } \rho = d, \\ O(n^{-\rho/d}), & \text{if } \rho < d, \end{cases}$$

where we can take $L(n) = \log n$, and omitting some straightforward details, we find that the assumed bounds for $\sigma_{AA}^{(n)}$ and $\sigma_{AB}^{(n)}$ imply that

$$n^{-2} \sum_{i,j\in\mathcal{I}_n} \sigma_{AA}^{(n)}(i-j)^r \sigma_{AB}^{(n)}(i-j)^{2(a_0-r)}$$

$$= \begin{cases} O(n^{-1}), & \text{if } p=0, (2-\alpha) > d, \\ O(n^{-1}\log n), & \text{if } p=0, (2-\alpha) = d, \\ O(n^{(\alpha-2)/d}), & \text{if } p=0, (2-\alpha) < d, \\ O(n^{-1}), & \text{if } p \geq 1, \end{cases}$$

where $d = 1, 2$. Note that the above statement is valid for each integer $m \geq 1$ and each integer $r$ satisfying $\max(0, a_0 - m) \leq r \leq a_0$. Also, using Lemma 4.2, it is straightforward to check that

$$\sum_{m=0}^{\infty} b_m^2 \left[\sum_{r=\max(a_0-m,0)}^{a_0} \{r!(a_0-r)!(a_0-r)!(m-a_0+r)!\}^{-1}\right] < \infty.$$

Finally, we put these results together and the proof is complete. □

The following result is elementary but is used repeatedly, and so is stated explicitly for convenience.

LEMMA 4.4. *Let $\{X_{k,n}: k = 1, \ldots, n; n \geq 1\}$ be an arbitrary triangular array of random variables such that $\sup_{k=1,\ldots,n} \sup_{n \geq 1} E|X_{k,n}| \leq C < \infty$. Then $n^{-1} \sum_{k=1}^{n} X_{k,n} = O_p(1)$ as $n \to \infty$.*

PROOF. Note that $E|n^{-1} \sum_k X_{k,n}| < C$ and then use the Markov inequality. □

The Prohorov metric $\rho$ and Ky Fan metric, here denoted $\kappa$, are defined as follows [see Dudley (1989)]. Let $X$ and $Y$ be random elements of a metric space $(S, \text{dist})$, with laws $P$ and $Q$, respectively, defined on the Borel sigma field of $(S, \text{dist})$. Then $\rho(P, Q) = \inf\{\varepsilon > 0 : P(A) \leq Q(A^\varepsilon) + \varepsilon \text{ for all Borelsets } A\}$, where $A^\varepsilon = \{y \in S : \text{dist}(x, y) < \varepsilon \text{ for some } x \in A\}$ and $\kappa(X, Y) = \inf\{\varepsilon > 0 : P[\text{dist}(X, Y) > \varepsilon] \leq \varepsilon\}$. Note that, by Theorem 11.3.5 of Dudley (1989), we have

$$\rho(P, Q) \leq \kappa(X, Y). \tag{4.3}$$

LEMMA 4.5. *Let $X$ and $Y$ be random elements with laws $P$ and $Q$, respectively, defined on the Borel sigma field of a metric space $(S, \text{dist})$, and let $\rho$ denote the Prohorov metric. Then $\rho(P, Q) \leq \{E[\text{dist}(X, Y)^2]\}^{2/3}$.*

PROOF. Chebyshev's inequality yields

$$P\Big[\text{dist}(X, Y) > \{E[\text{dist}(X, Y)^2]\}^{1/3}\Big] \leq \{E[\text{dist}(X, Y)^2]\}^{1/3},$$

and so the result follows from (4.3). □

In all applications of this result given below, $S = \mathbf{R}^2$ and dist is the usual Euclidean metric.

**5. Proofs.** We now prove the theorems stated in Section 2. The structure of the proof of each theorem is very similar, and, in fact, the proofs given for Steps 1–5 cover both theorems. The only substantial difference between Theorems A and B is in the limit distribution which arises in Step 6. Throughout the proof, we will use the multi-index notation specified in the Appendix, on the understanding that $d = 1$ or 2.

We first introduce some notation that will be used throughout this section. Write

$$W_{ij} = n^{\alpha/(2d)} \left\{ X\left(\frac{i+j}{n_0}\right) - X\left(\frac{i}{n_0}\right) \right\}$$

and

$$(5.1) \quad Y_{iu} = n^{\alpha/(2d)} \sum_j a_j^u X\left(\frac{i+j}{n_0}\right) = \sum_j a_j^u W_{ij},$$

where $X(t)$ is the underlying Gaussian field (see Section 2). Note that the last equality is a consequence of the fact that $\sum_j a_j^u = 0$ (see the Appendix). Define

$$\sigma_{WW}(i-j) = \sup_{k,l} |\text{cov}(W_{ik}, W_{jl})|,$$

$$\sigma_{WX}(i-j) = \sup_k |\text{cov}\{W_{ik}, X(j/n_0)\}|,$$

$$\sigma_{YY}(i-j) = \sup_{u,v} |\text{cov}(Y_{iu}, Y_{jv})|,$$

$$\sigma_{YX}(i-j) = \sup_u |\text{cov}\{Y_{iu}, X(j/n_0)\}|.$$

Then condition $(\mathcal{A}1)_4^{(d)}$ implies the existence of a constant $C$ independent of $i$, $j$ and $n$ such that

$$(5.2) \quad \begin{aligned} \sigma_{WW}(i-j) &\leq C(1+|i-j|)^{\alpha-2}, \\ \sigma_{WX}(i-j) &\leq Cn^{-\alpha/(2d)}(1+|i-j|)^{\alpha-1}, \end{aligned}$$

$$(5.3) \quad \begin{aligned} \sigma_{YY}(i-j) &\leq C(1+|i-j|)^{\alpha-2p-2}, \\ \sigma_{YX}(i-j) &\leq Cn^{-\alpha/(2d)}(1+|i-j|)^{\alpha-p-1}, \end{aligned}$$

where $p$ is the order of the increment $\mathbf{a}$ on which $Y$ is based. See Kent and Wood (1997) and Chan and Wood (2000), Lemma 3.1, for justification of (5.2) and (5.3).

PROOF OF THEOREM A. The proof is broken into a number of steps. The $T_i$ and $T_{ij}$ referred to below are defined in the course of the proof; each of these quantities is $O_p(1)$ and in some cases of smaller order.

STEP 1. Show that it is sufficient to prove the theorem for those $g$ which satisfy $(\mathcal{A}2)$ and $(\mathcal{A}3)$ and have compact support.

STEP 2. Show that

$$\bar{G}_1(\hat{\alpha} - \alpha_n) = T_0 + n^{-\alpha/(2d)}T_1 + n^{-\alpha/d}T_2 + O_p(n^{-3\alpha/(2d)} + n^{-1}),$$

where

$$(5.4) \quad \bar{G}_1 = n^{-1} \sum_{i \in \mathcal{I}_n} [g^{(1)}\{X(i/n_0)\}]^2$$

and $T_0$, $T_1$ and $T_2$ are defined in (5.9) via (5.6)–(5.8).

STEP 3. Show that

$$T_0 = T_{00} + n^{-\alpha/d}T_{01} + O_p(n^{-3\alpha/(2d)})$$
$$+ \begin{cases} O(n^{-1/2}), & \text{if } p = 0 \text{ and } d = 2, \\ O(n^{-1}), & \text{if } p \geq 1 \text{ and/or } d = 1, \end{cases}$$

where $T_{00}$ and $T_{01}$ are defined in (5.16) and (5.14), respectively, and show that

$$T_{01} \xrightarrow{\mathcal{D}} \left( \sum_{u=1}^{m} L_u \mu_{0,u}^{-1} \tau_{0,1u} \right)$$
$$\times \int_{t \in [0,1]^d} [g^{(1)}\{X(t)\}g^{(3)}\{X(t)\} + [g^{(2)}\{X(t)\}]^2] \, dt$$

as $n \to \infty$, where $\mu_{0,u}$ and $\tau_{0,1u}$ are defined below (2.9).

STEP 4. Show that $T_1 = n^{-\alpha/(2d)}T_{11} + O_p(n^{-\alpha/d})$, where $T_{11}$ is defined in (5.17), and show that as $n \to \infty$,

$$T_{11} \xrightarrow{\mathcal{D}} \left( \sum_{u=1}^{m} L_u \mu_{0,u}^{-1} \tau_{0,2u} \right)$$
$$\times \int_{t \in [0,1]^d} [g^{(1)}\{X(t)\}g^{(3)}\{X(t)\} + [g^{(2)}\{X(t)\}]^2] \, dt,$$

where $\tau_{0,2u}$ is defined below (2.9).

STEP 5. Show that $T_2 = T_{21} + T_{22} + O_p(n^{-\alpha/(2d)} + n^{-1/2})$, where $T_{21}$ and $T_{22}$ are defined in (5.19) and (5.20), respectively, and show that, as $n \to \infty$,

$$T_{21} \xrightarrow{\mathcal{D}} \left( \sum_{u=1}^{m} L_u \mu_{0,u}^{-1} \tau_{0,3u} \right) \int_{t \in [0,1]^d} g^{(1)}\{X(t)\}g^{(3)}\{X(t)\} \, dt$$

and

$$T_{22} \xrightarrow{\mathcal{D}} \left( \sum_{u=1}^{m} L_u \mu_{0,u}^{-1} \tau_{0,4u} \right) \int_{t \in [0,1]^d} [g^{(2)}\{X(t)\}]^2 \, dt,$$

and $\tau_{0,3u}$ and $\tau_{0,4u}$ are defined below (2.9).

STEP 6. Establish convergence in distribution of $(n^{1/2}T_{00}, \bar{G}_1)$, where $T_{00}$ is defined in (5.16) and $\bar{G}_1$ is defined in (5.4).

PROOF OF STEP 1. Condition $(\mathcal{A}1)_4^{(d)}$ implies (1.1) when $d = 1$ and (2.3) when $d = 2$. In each case, we may use Kolmogorov's lemma [see, e.g., Rogers and Williams (1994), page 59] to establish that $X(t)$ has a continuous version on $[0,1]^d$. Consequently, for each $\varepsilon > 0$ there exists a $C$, depending on $\varepsilon$ and the distribution of $X$, such that

$$(5.5) \qquad P\left[\sup_{t \in [0,1]^d} |X(t)| > C\right] < \varepsilon.$$

For given $g$, let $g_C$ denote a function with compact support such that $g(t) = g_C(t)$ for all $\|t\| < C$, and let $\hat{\alpha}_C$ denote the estimator of $\alpha$ that would have been obtained if $g_C\{X(t)\}$ rather than $g\{X(t)\}$ had been observed. It follows from (5.5) that

$$P[\hat{\alpha} \neq \hat{\alpha}_C \text{ for some } n] < \varepsilon.$$

As a consequence, if the theorem is true for all functions of compact support which satisfy assumptions $(\mathcal{A}2)$ and $(\mathcal{A}3)$, then it is also true for each $g$ which satisfies $(\mathcal{A}2)$ and $(\mathcal{A}3)$, whether or not $g$ has compact support. This argument can be established rigorously using probability metrics (cf. the argument given in Step 6). We omit the details. □

For the remainder of the proof we shall assume that $g$ has compact support [in addition to satisfying $(\mathcal{A}2)$ and $(\mathcal{A}3)$].

PROOF OF STEP 2. By Taylor's theorem

$$\sum_j a_j^u g_{i+j} = \sum_j a_j^u (g_{i+j} - g_i) \qquad \left(\text{since } \sum_j a_j^u = 0\right)$$

$$= \sum_j a_j^u \left(\sum_{r=1}^3 n^{-r\alpha/(2d)} (r!)^{-1} W_{ij}^r g_i^{(r)} + n^{-2\alpha/d}(4!)^{-1} W_{ij}^4 \tilde{g}_i^{(4)}\right)$$

$$= \sum_{r=1}^4 n^{-r\alpha/(2d)} M_{riu},$$

where, for $r = 1, 2, 3$,

$$M_{riu} = g_i^{(r)} \frac{1}{r!} \sum_j a_j^u W_{ij}^r, \qquad M_{4iu} = (4!)^{-1} \sum_j a_j^u \tilde{g}_{ij}^{(4)} W_{ij}^4,$$

$g_i^{(r)}$ is $g^{(r)}$ evaluated at $X(i/n_0)$, $r = 1, 2, 3$, and from Taylor's theorem, $\tilde{g}_{ij}^{(4)} = g_i^{(4)}[\theta_j X\{(i+j)/n_0\} + (1-\theta_j)X(i/n_0)]$, where each $\theta_j \in [0,1]$ is suitably chosen.

Then
$$n^{\alpha/d}\bar{Z}_u = n^{-1}\sum_{i\in\mathcal{I}_n}\left(n^{\alpha/(2d)}\sum_j a_j^u g_{i+j}\right)^2$$
$$= n^{-1}\sum_{i\in\mathcal{I}_n}\left(\sum_{r=1}^4 n^{-(r-1)\alpha/(2d)}M_{riu}\right)^2$$
$$= \left(n^{-1}\sum_{i\in\mathcal{I}_n}M_{1iu}^2\right) + n^{-\alpha/(2d)}\mu_u\tilde{T}_{1u} + n^{-\alpha/d}\mu_u\tilde{T}_{2u}$$
$$+ n^{-3\alpha/(2d)}n^{-1}\sum_{i\in\mathcal{I}_n}S_{iu},$$

where

(5.6) $$\tilde{T}_{1u} = 2\mu_u^{-1}n^{-1}\sum_{i\in\mathcal{I}_n}M_{1iu}M_{2iu},$$

(5.7) $$\tilde{T}_{2u} = \mu_u^{-1}n^{-1}\sum_{i\in\mathcal{I}_n}(M_{2iu}^2 + 2M_{1iu}M_{3iu})$$

and
$$S_{iu} = 2(M_{1iu}M_{4iu} + M_{2iu}M_{3iu}) + n^{-\alpha/(2d)}(M_{3iu}^2 + 2M_{2iu}M_{4iu})$$
$$+ 2n^{-\alpha/d}M_{3iu} + M_{4iu} + n^{-3\alpha/(2d)}M_{4iu}^2.$$

Since by Step 1 we are assuming that the $g^{(r)}$ are bounded for $1\le r\le 4$, it is a straightforward (if tedious) matter to check that $E|S_{iu}| < \infty$. Therefore, since the $S_{iu}$ are identically distributed for each $1\le u\le m$, it follows from Lemma 4.4 that
$$n^{-3\alpha/(2d)}n^{-1}\sum_{i\in\mathcal{I}_n}S_{iu} = O_p(n^{-3\alpha/(2d)})$$

for each $u$. Moreover,
$$n^{-1}\sum_{i\in\mathcal{I}_n}M_{1iu}^2 = n^{-1}\sum_{i\in\mathcal{I}_n}\left(\sum_j a_j^u W_{ij}\right)^2 (g_i^{(1)})^2$$
$$= n^{-1}\sum_{i\in\mathcal{I}_n}(Y_{iu}^2 - \mu_u + \mu_u)(g_i^{(1)})^2$$
$$= \mu_u\bar{G}_1 + \mu_u\tilde{T}_{0u},$$

where $\bar{G}_1$ is defined in (5.4) and

(5.8) $$\tilde{T}_{0u} = n^{-1}\sum_{i\in\mathcal{I}_n}(\mu_u^{-1}Y_{iu}^2 - 1)\{g_i^{(1)}\}^2.$$

So $n^{\alpha/d}\bar{Z}_u = \bar{G}_1\mu_u + R_u\mu_u$, where

$$R_u = \tilde{T}_{0u} + n^{-\alpha/(2d)}\tilde{T}_{1u} + n^{-\alpha/d}\tilde{T}_{2u} + O_p(n^{-3\alpha/(2d)}).$$

Since each $\tilde{T}_{ku}$ is bounded in probability for $k=1,2,3$, it follows that

$$\begin{aligned}
\bar{G}_1(\hat{\alpha} - \alpha_n) &= \bar{G}_1 \sum_{u=1}^{m} L_u \log(\bar{Z}_u/\mu_u) \\
&= \bar{G}_1 \sum_{u=1}^{m} L_u\{\log \bar{G}_1 + \log(1 + R_u/\bar{G}_1)\} \\
&= \sum_{u=1}^{m} L_u R_u\{1 + O_p(R_u/\bar{G}_1)\} \\
&= \sum_{u=1}^{m} L_u[\tilde{T}_{0u} + n^{-\alpha/(2d)}\tilde{T}_{1u} + n^{-\alpha/d}\tilde{T}_{2u} + O_p\{n^{-3\alpha/(2d)} + R_u^2\}] \\
&= T_0 + n^{-\alpha/(2d)} T_1 + n^{-\alpha/d} T_2 + O_p(n^{-3\alpha/(2d)} + n^{-1}),
\end{aligned}$$

where

(5.9) $$T_k = \sum_{u=1}^{m} L_u \tilde{T}_{ku}, \qquad k=0,1,2,$$

and we have used the fact that $R_u = O_p(n^{-1/2} + n^{-\alpha/d})$, so that $R_u^2 = O_p(n^{-1} + n^{-2\alpha/d})$. The order statement for $R_u$ follows from Steps 3–6. □

We now introduce some notation which is needed in Steps 3–5. Recall the definition of $W_{ij}$ given at the beginning of Section 5. Writing $\zeta_j^{(n)} = \text{cov}(X(i/n_0), W_{ij})$, define $\zeta^{(n)} = (\zeta_j^{(n)}, -mJ \leq j \leq mJ)$ and $W_i = (W_{ij}, -mJ \leq j \leq mJ)$, and let $V_W^{(n)}$ denote the covariance matrix of $W_i$. Note that $V_W^{(n)}$ does not depend on $i$ because of the stationarity of $W_i$, but that the distribution of $W_i$ does depend on $n$; this dependence on $n$ has been suppressed for notational convenience. Define $b^{(n)} = (b_j^{(n)}, -mJ \leq j \leq mJ)$ by $b^{(n)} = n^{\alpha/(2d)}(V_W^{(n)})^{-1}\zeta^{(n)}$. Note that $n^{-\alpha/(2d)}\sum_j b_j^{(n)} W_{ij}$ is the projection of $X(i/n)$ onto the span of $W_{ij}, -mJ \leq j \leq mJ$. Let $b^{(0)} = \lim_{n\to\infty} b^{(n)}$ denote the limit of $b^{(n)}$ which necessarily exists under assumption $(\mathcal{A}1)_4^{(d)}$.

We may write

$$X(i/n_0) = n^{-\alpha/(2d)}\left(\sum_j b_j^{(n)} W_{ij}\right) + (1 - n^{-\alpha/d} b^{(n)T} V_W^{(n)} b^{(n)}/\gamma_0)^{1/2} \check{X}_i,$$

where $\check{X}_i \equiv \check{X}(i/n_0)$, $\check{X}_i \sim N(0, \gamma_0)$ is independent of $W_i$, $\gamma_0 = \gamma(0)$ is the variance of $X(i/n_0)$ and $V_W^{(n)}$ and $b^{(n)}$ are as defined above. Then

$$X(i/n_0) = \check{X}_i + n^{-\alpha/(2d)} \delta_{1i} + n^{-\alpha/d} \delta_{2i} + n^{-2\alpha/d} \delta_{3i}$$
$$= \check{X}_i + n^{-\alpha/(2d)} \delta_{0i},$$

where

$$\delta_{1i} = \sum_j b_j^{(n)} W_{ij}, \qquad \delta_{2i} = -\{b^{(n)T} V_W^{(n)} b^{(n)} / (2\gamma_0)\} \check{X}_i,$$

$$n^{-2\alpha/d} \delta_{3i} = (1 - n^{-\alpha/d} b^{(n)T} V_W^{(n)} b^{(n)} / \gamma_0)^{1/2} \check{X}_i - \check{X}_i - n^{-\alpha/d} \delta_{2i}$$

and

$$\delta_{0i} = \delta_{1i} + n^{-\alpha/(2d)} \delta_{2i} + n^{-3\alpha/(2d)} \delta_{3i}.$$

PROOF OF STEP 3. Writing $g_i^{(r)}$ for $g^{(r)}(X(i/n_0))$ as before, and $\check{g}_i^{(r)}$ for $g^{(r)}(\check{X}_i)$, we obtain

$$(5.10) \quad \begin{aligned} (g_i^{(1)})^2 &= (\check{g}_i^{(1)})^2 + 2n^{-\alpha/(2d)} \delta_{0i} \check{g}_i^{(1)} \check{g}_i^{(2)} \\ &\quad + n^{-\alpha/d} \delta_{0i}^2 \{\check{g}_i^{(1)} \check{g}_i^{(3)} + (\check{g}_i^{(2)})^2\} + n^{-3\alpha/(2d)} \delta_{0i}^3 R_{1i}, \end{aligned}$$

where $R_{1i}$ is a remainder term which can be determined explicitly.

We now study the contribution of each of the four terms on the right-hand side of (5.10). First, note that $\delta_{0i}^3 R_{1i}$ can be expressed as a finite sum of terms, each of which can be expressed as a bounded function multiplied by a polynomial in Gaussian variables. Therefore, since the $\delta_{0i}^3 R_{1i}$ ($i \in \mathcal{I}_n$) are identically distributed, we may use Lemma 4.4 to show that

$$(5.11) \quad \sum_{u=1}^m L_u \mu_u^{-1} \left( n^{-1} \sum_{i \in \mathcal{I}_n} (n^{-3\alpha/(2d)} \delta_{0i}^3 R_{1i}) \right) = O_p(n^{-3\alpha/(2d)}).$$

Also, using similar arguments,

$$(5.12) \quad \begin{aligned} \left(\sum_j a_j^u W_{ij}\right)^2 \delta_{0i}^2 \{\check{g}_i^{(1)} \check{g}_i^{(3)} + (\check{g}_i^{(2)})^2\} \\ = \tau_{1u} \{\check{g}_i^{(1)} \check{g}_i^{(3)} + (\check{g}_i^{(2)})^2\} \\ + \left\{\left(\sum_j a_j^u W_{ij}\right)^2 \delta_{0i}^2 - \tau_{1u}\right\} \{\check{g}_i^{(1)} \check{g}_i^{(3)} + (\check{g}_i^{(2)})^2\}, \end{aligned}$$

where

$$(5.13) \quad \tau_{1u} = E\left[\left(\sum_j a_j^u W_{ij}\right)^2 \left(\sum_k b_k^{(n)} W_{ik}\right)^2\right].$$

It follows, after applying Lemmas 4.3 and 4.4 to the second term on the right-hand side of (5.12), that

$$\sum_{u=1}^{m} L_u \mu_u^{-1} \left( n^{-1} \sum_{i \in \mathcal{I}_n} \left( \sum_j a_j^u W_{ij} \right)^2 \delta_{0i}^2 \{\breve{g}_i^{(1)} \breve{g}_i^{(3)} + (\breve{g}_i^{(2)})^2\} \right)$$
$$= T_{01} + O_p(n^{-\alpha/(2d)}),$$

where

$$(5.14) \quad T_{01} = \left( \sum_{u=1}^{m} L_u \mu_u^{-1} \tau_{1u} \right) \left( n^{-1} \sum_{i \in \mathcal{I}_n} \breve{g}_i^{(1)} \breve{g}_i^{(3)} + (\breve{g}_i^{(2)})^2 \right).$$

A similar argument, using Lemmas 4.3 and 4.4 again, shows that

$$(5.15) \quad n^{-\alpha/(2d)} n^{-1} \sum_{i \in \mathcal{I}_n} \left( \sum_j a_j^u W_{ij} \right)^2 \delta_{0i} \breve{g}_i^{(1)} \breve{g}_i^{(2)} = O_p(n^{-1}),$$

except when $p = 0$ and $d = 2$, in which case we can conclude that the left-hand side of (5.15) is of size $O_p(n^{-1/2})$. Writing

$$(5.16) \quad T_{00} = \sum_{u=1}^{m} L_u n^{-1} \sum_{i \in \mathcal{I}_n} (\mu_u^{-1} Y_{iu}^2 - 1)(\breve{g}_i^{(1)})^2$$

and putting (5.10)–(5.16) together, we see that the first part of Step 3 is proved. To establish the limiting distribution of $T_{01}$, we may use an identical argument to that given in Step 6 to prove that $J_n(h) \xrightarrow{\mathcal{D}} J_0(h)$ [which is the more straightforward part of showing that $K_n(h) \xrightarrow{\mathcal{D}} K_0(h)$]. To avoid duplication, we omit the details. □

PROOF OF STEP 4. By a Taylor expansion,

$$g_i^{(1)} g_i^{(2)} = \breve{g}_i^{(1)} \breve{g}_i^{(2)} + n^{-\alpha/(2d)} \delta_{0i} \{\breve{g}_i^{(1)} \breve{g}_i^{(3)} + (\breve{g}_i^{(2)})^2\} + n^{-\alpha/d} \delta_{0i}^2 R_{3i}.$$

Using similar arguments to those in Step 3, in particular Lemma 4.3, we find that for any $\varepsilon > 0$,

$$T_1 = \sum_{u=1}^{m} L_u \mu_u^{-1} \left( n^{-1} \sum_{i \in \mathcal{I}_n} \left( \sum_j a_j^u W_{ij} \right) \left( \sum_k a_k^u W_{ik}^2 \right) g_i^{(1)} g_i^{(2)} \right)$$
$$= n^{-\alpha/(2d)} T_{11} + O_p(n^{-\alpha/d}) + O(n^{-1+\varepsilon}),$$

where

$$(5.17) \quad T_{11} = \left( \sum_{u=1}^{m} L_u \mu_u^{-1} \tau_{2u} \right) n^{-1} \sum_{i \in \mathcal{I}_n} \{\breve{g}_i^{(1)} \breve{g}_i^{(3)} + (\breve{g}_i^{(2)})^2\}$$

and

$$\tau_{2u} = E\left\{\left(\sum_j a_j^u W_{ij}\right)\left(\sum_k a_k^u W_{ik}^2\right)\left(\sum_l b_l^{(n)} W_{il}\right)\right\}, \tag{5.18}$$

so the first part of Step 4 is proved. Justification of the claim concerning the limit distribution of $T_{11}$ follows along similar lines to that for $T_{01}$ in Step 3. □

PROOF OF STEP 5. In this case
$g_i^{(1)} g_i^{(3)} = \breve{g}_i^{(1)} \breve{g}_i^{(3)} + n^{-\alpha/(2d)} \delta_{0i} R_{4i}$ and $(g_i^{(2)})^2 = (\breve{g}_i^{(2)})^2 + n^{-\alpha/(2d)} \delta_{0i} R_{5i}$.
Using Lemmas 4.3 and 4.4 again and writing

$$T_{21} = \left(\sum_{u=1}^m L_u \mu_u^{-1} \tau_{3u}\right) n^{-1} \sum_{i \in \mathcal{I}_n} \breve{g}_i^{(1)} \breve{g}_i^{(3)}, \tag{5.19}$$

$$T_{22} = \left(\sum_{u=1}^m L_u \mu_u^{-1} \tau_{4u}\right) n^{-1} \sum_{i \in \mathcal{I}_n} (\breve{g}_i^{(2)})^2, \tag{5.20}$$

$$\tau_{3u} = \tfrac{1}{3} E\left(\sum_j a_j^u W_{ij}\right)^2 \left(\sum_j a_j^u W_{ij}^3\right), \qquad \tau_{4u} = \tfrac{1}{4} E\left(\sum_j a_j^u W_{ij}^2\right)^2, \tag{5.21}$$

we find that $n^{-\alpha/d} T_2 = n^{-\alpha/d}(T_{21} + T_{22}) + O_p(n^{-1})$ as required. Justification of the claim concerning the limit distributions of $T_{21}$ and $T_{22}$ follows along similar lines to that for $T_{01}$ in Step 3 and $T_{11}$ in Step 4. □

PROOF OF STEP 6. Here we shall show that

$$(n^{1/2} T_{00}, \bar{G}_1) \xrightarrow{\mathcal{D}} \left(\sigma \int_{[0,1]^d} [g^{(1)}\{X(t)\}]^2 \, dB(t), G_1\right), \tag{5.22}$$

where $T_{00}$ is given in (5.16), $\bar{G}_1$ is defined in (5.4), $G_1$ is defined in (2.5) with $r=1$, $\{B(t)\}$ is the random Gaussian measure given in (2.11), which is independent of $\{X(t)\}$, and $\sigma$ is given in (2.10). Recall that $n_0 = (n_0[1], \ldots, n_0[d])$ and $n = \prod_{l=1}^d n_0[l]$ and that we assume that condition $(\mathcal{A}4)$ in Section 2.1 holds when $d = 2$.

Let $\mathcal{H}$ denote the class of smooth functions with compact support. For $h \in \mathcal{H}$, define

$$I_n(h) = n^{-1/2} \sum_{i \in \mathcal{I}_n} \xi_i h\{\breve{X}(i/n_0)\}, \qquad J_n(h) = n^{-1} \sum_{i \in \mathcal{I}_n} h\{X(i/n_0)\}$$

and $K_n(h) = (I_n(h), J_n(h))$, where $\xi_i = \sum_{u=1}^m L_u(\mu_u^{-1} Y_{iu}^2 - 1)$. Note that by construction $K_n\{(g^{(1)})^2\} = (n^{1/2} T_{00}, \bar{G}_1)$. Also define

$$I_0(h) = \int_{[0,1]^d} h\{X(t)\} \, dB(t), \qquad J_0(h) = \int_{[0,1]^d} h\{X(t)\} \, dt$$

and $K_0(h) = (I_0(h), J_0(h))$. Note that by construction $K_0\{(g^{(1)})^2\}$ is equal to the right-hand side of (5.22).

We will show that, for each $h \in \mathcal{H}$, $K_n(h)$ converges to $K_0(h)$ in distribution. Then, in view of (5.22) and Step 1, Step 6 will follow.

For each positive integer $r$, let $\pi_r$ denote a partition of $[0,1]^d$ given by $\pi_r = \{A_j, j \in \mathcal{J}_r\}$, where $\mathcal{J}_r = \{j \in \mathbf{Z}^d : 0 \le j[l] < r, 1 \le l \le d\}$, where $A_j \subset [0,1]^d$ is defined in the following way. Write $t_j = r^{-1}j$. When $d = 1$, $A_j$ is the interval of width $r^{-1}$ centered at $t_j = r^{-1}j + \frac{1}{2}$; when $d = 2$, $A_j$ is a square with sides of length $r^{-1}$ which are parallel to the coordinate axes, and each $A_j$ is centered at $t_j = r^{-1}j + \frac{1}{2}\mathbf{1}$, where $\mathbf{1} = (1,1)$ is a 2-vector of 1's. Given $\pi_r$, we define two functions, $t^*(t)$ and $i^*(i)$, as follows: for $t \in [0,1]^d$, $t^*(t) = t_j$ when $t \in A_j$, and for each multi-index $0 \le i < n_0$, we define $i^*(i) \equiv i_n^*(i) = n_0 t_j$, where $j$ is such that $i/n_0 \in A_j$. Define

$$I_n^*(h) \equiv I_n^*(h; \pi_r) = n^{-1/2} \sum_{i \in \mathcal{I}_n} \xi_i h\{\check{X}(i^*/n_0)\},$$

$$J_n^*(h) \equiv J_n^*(h; \pi_r) = n^{-1} \sum_{i \in \mathcal{I}_n} h\{X(i^*/n_0)\}$$

and $K_n^*(h) = (I_n^*(h), J_n^*(h))$, and write

$$I_0^*(h) \equiv I_0^*(h; \pi_r) = \sigma \int_{[0,1]^d} h\{X(t^*)\} \, dB(t),$$

$$J_0^* \equiv J_0^*(h; \pi_r) = \int_{[0,1]^d} h\{X(t^*)\} \, dt$$

and $K_0^*(h) = (I_0^*(h), J_0^*(h))$, where $\{B(t)\}$ is a random Gaussian measure of the form (2.11) which is independent of $\{X(t)\}$.

Let $P_{n,h}$, $P_{0,h}$, $P_{n,h}^*$ and $P_{0,h}^*$ denote the distributions of $K_n(h)$, $K_0(h)$, $K_n^*(h)$ and $K_0^*(h)$, respectively. We need to show that $P_{n,h} \xrightarrow{\mathcal{D}} P_{0,h}$. We shall do this by showing that, given $\varepsilon > 0$, there exists a partition $\pi_r$ of $[0,1]^d$, depending on $\varepsilon$ and $h$, such that

(5.23) $$\rho(P_{0,h}, P_{0,h}^*) < \varepsilon/3$$

and, for any such $\pi_r$, when $n$ is sufficiently large,

(5.24) $$\rho(P_{n,h}^*, P_{0,h}^*) < \varepsilon/3;$$

and for a sufficiently fine partition $\pi_r$ and $n$ sufficiently large,

(5.25) $$\rho(P_{n,h}, P_{n,h}^*) < \varepsilon/3,$$

where $\rho$ denotes the Prohorov metric. Then, when $n$ is sufficiently large,

$$\rho(P_{n,h}, P_{0,h}) \le \rho(P_{n,h}, P_{n,h}^*) + \rho(P_{n,h}^*, P_{0,h}^*) + \rho(P_{0,h}^*, P_{0,h}) < \varepsilon,$$

and, since $\varepsilon > 0$ may be chosen arbitrarily small, Step 6 will have been proved.

PROOF OF (5.23). Using Lemma 4.5,

$$\rho(P_{0,h}, P_{0,h}^*) \leq [E\{I_0(h) - I_0^*(h)\}^2 + E\{J_0(h) - J_0^*(h)\}^2]^{2/3}. \quad (5.26)$$

Using Fubini's theorem, the Cauchy–Schwarz inequality and the fact that $\{X(t)\}$ and $\{B(t)\}$ are independent, we obtain

$$E\{I_0(h) - I_0^*(h)\}^2 = E\left[\sigma \int_{[0,1]^d} [h\{X(t)\} - h\{X(t^*)\}]\, dB(t)\right]^2$$
$$= \sigma^2 \int_{[0,1]^d} E[h\{X(t)\} - h\{X(t^*)\}]^2\, dt \quad (5.27)$$

and

$$E\{J_0(h) - J_0^*(h)\}^2 \leq \int_{[0,1]^d} E[h\{X(t)\} - h\{X(t^*)\}]^2\, dt. \quad (5.28)$$

Given $h \in \mathcal{H}$ and $\varepsilon > 0$, it is clear (since $\{X(t)\}$ is almost surely continuous and therefore uniformly continuous on $[0,1]^d$) that we can find a (sufficiently fine) partition $\pi_r$ of $[0,1]^d$ such that

$$\int_{[0,1]^d} E[h\{X(t)\} - h\{X(t^*)\}]^2 < \tfrac{1}{2}(\varepsilon/3)^{3/2}.$$

Then, using (5.27) and (5.28),

$$E\{I_0(h) - I_0^*(h)\}^2 < \tfrac{1}{2}(\varepsilon/3)^{3/2} \quad \text{and} \quad E\{J_0(h) - J_0^*(h)\}^2 < \tfrac{1}{2}(\varepsilon/3)^{3/2},$$

in which case (5.23) follows from (5.26). □

PROOF OF (5.24). For any Lebesgue-measurable set $A \subseteq [0,1]^d$, define

$$\Xi_n(A) = n^{-1/2} \sum_{i: i/n_0 \in A} \xi_i. \quad (5.29)$$

Then, for a given partition $\pi_r$,

$$I_n^*(h) \equiv n^{-1/2} \sum_{i \in \mathcal{I}_n} \xi_i h\{\check{X}(i^*/n_0)\} = \sum_{j \in \mathcal{J}_r} \Xi_n(A_j) h\{\check{X}(t_j)\}$$

and

$$J_n^*(h) \equiv n^{-1} \sum_{i \in \mathcal{I}_n} h\{X(i^*/n_0)\} = n^{-1} \sum_{j \in \mathcal{J}_r} c_j h\{X(t_j)\},$$

where $c_j = \#\{i : i/n_0 \in A_j\} \sim n/r^d$. In all cases covered by Theorem A a central limit theorem applies to $\Xi_n([0,1]^d)$; that is, $\Xi_n([0,1]^d) \xrightarrow{\mathcal{D}} N(0, \sigma^2)$, where $\sigma^2$ is defined in (2.10). See Kent and Wood (1995, 1997) and Chan and Wood (2000) for further details. It follows from a straightforward extension of those proofs that, for any fixed $\pi_r$, $\{\Xi_n(A_j), j \in \mathcal{J}_r\} \xrightarrow{\mathcal{D}} \{\sigma B(A_j), j \in \mathcal{J}_r\}$, where $B$ is the random Gaussian measure defined in (2.11). Moreover, it is an easy consequence of the definition of $\check{X}$ that $\{\check{X}(t_j), j \in \mathcal{J}_r\} \xrightarrow{\mathcal{D}} \{X(t_j), j \in \mathcal{J}_r\}$, where the random variables $\{X(t_j), j \in \mathcal{J}_r\}$ are independent of the random variables $\{B(A_j), j \in \mathcal{J}_r\}$. Consequently, an application of the continuous mapping theorem implies that $K_n^*(h) \xrightarrow{\mathcal{D}} K_0^*(h)$, from which (5.24) follows for sufficiently large $n$, since the Prohorov metric $\rho$ metrizes convergence in distribution. □

PROOF OF (5.25). We will only sketch the proof of this result. It is sufficient to show that

(5.30) $$\lim_{n \to \infty} E\{I_n(h) - I_n^*(h)\}^2 = \sigma^2 \int_{[0,1]^d} E[h\{X(t)\} - h\{X(t^*)\}]^2 \, dt$$

and

(5.31) $$\lim_{n \to \infty} E\{J_n(h) - J_n^*(h)\}^2 \leq \int_{[0,1]^d} E[h\{X(t)\} - h\{X(t^*)\}]^2 \, dt,$$

because, if we choose a partition $\pi$ such that, for all $n$ sufficiently large,

$$\int_{[0,1]^d} E[h\{X(t)\} - h\{X(t^*)\}]^2 \, dt < \tfrac{1}{2}(\varepsilon/6)^{3/2},$$

then (5.25) will follow from (5.30), (5.31) and Lemma 4.5.

The proof of (5.31) is omitted as it is quite straightforward. However, we outline the proof of (5.30), which requires more work. For simplicity, we focus on the case in which $h = \sum_m c_m H_m$ is a polynomial; it is an easy matter to extend the proof to all $h \in \mathcal{H}$. Now

$$E\{I_n(h) - I_n^*(h)\}^2 = n^{-1} \sum_{i,j \in \mathcal{I}_n} E[\xi_i \xi_j (h_i - h_{i^*})(h_j - h_{j^*})]$$

$$= \sum_{m,m'=1}^{M} \frac{c_m c_{m'}}{m! m'!} n^{-1} \sum_{i,j \in \mathcal{I}_n} \delta(i, j, m, m'),$$

where $h_i = h\{\check{X}(i/n_0)\}$, $h_{i^*} = h\{\check{X}(i^*/n_0)\}$ and

$$\delta(i, j, m, m') = E[\xi_i \xi_j \{H_m(\check{X}_i) - H_m(\check{X}_{i^*})\} \{H_{m'}(\check{X}_j) - H_{m'}(\check{X}_{j^*})\}].$$

Then (5.30) is a consequence of the following: for each $m \neq m'$,

(5.32) $$n^{-1} \sum_{i,j \in \mathcal{I}_n} \delta(i, j, m, m') = o(1);$$

and

$$\lim_{n\to\infty} n^{-1} \sum_{i,j\in\mathcal{I}_n} \delta(i,j,m,m)$$
$$\times \lim_{n\to\infty} n^{-1} \sum_{i,j\in\mathcal{I}_n} E[\xi_i \xi_j] E[H_m(\check{X}_i) - H_m(\check{X}_{i^*})]$$
(5.33)
$$\times [H_m(\check{X}_j) - H_m(\check{X}_{j^*})]$$
$$= \sigma^2 \int_{[0,1]^d} [H_m(X(t)) - H_m(X(t^*))]^2 \, dt.$$

The results (5.32) and (5.33) are derived using the diagram formula (see the references given in Remark 4.1) combined with (5.2) and (5.3). The arguments are broadly similar to those used in the proof of Lemma 4.3. That concludes our sketched proof of (5.25). □

Thus, Step 6 is now complete. □

To conclude the proof of Theorem A, we use (5.22) and the continuous mapping theorem to show that

$$\bar{G}_1^{-1} n^{1/2} T_{00} \xrightarrow{\mathcal{D}} \frac{\sigma}{G_1} \int_{[0,1]^d} [g^{(1)}\{X(t)\}]^2 \, dB(t) \xrightarrow{\mathcal{D}} \sigma \frac{\sqrt{G_2}}{G_1} Z,$$

where $Z \sim N(0,1)$ is independent of $G_1$ and $G_2$, and $G_2$ is defined in (2.5) with $r = 2$. Finally, putting Steps 1–6 together, we see that the proof of Theorem A is complete. □

PROOF OF THEOREM B. The proof of Theorem B is essentially the same as the proof of Theorem A, except that Step 6 is different. The decomposition given by (5.23)–(5.25) can still be used, but the principal difference is that $\Xi_n$ in (5.29) does not obey a central limit theorem, and therefore $I_n(h)$ does not converge to a stochastic integral with respect to the random Gaussian measure $B(t)$. What actually happens is that, when $d = 1$ and $3/2 < \alpha < 2$ or $d = 2$ and $1 < \alpha < 2$,

$$\operatorname{var}\left(n^{-1} \sum_{i\in\mathcal{I}_n} \xi\right) = O(n^{(2\alpha-4)/d}),$$

and $n^{(2-\alpha)/d} n^{-1} \sum_{i\in\mathcal{I}_n} \xi_i$ obeys a noncentral limit theorem of the type given by Theorem 6 in Arcones (1994). Then it can be shown, via the decomposition (5.23)–(5.25), that $I_n(h)$ converges to the Wiener–Itô integral specified in the statement of Theorem B. The particular form of the spectral measure $S$ follows in part as a consequence of Theorem 1' of Dobrushin and Major (1979) and in part as a consequence of the degeneracy result given in part (c) of Theorem 1 in Kent and Wood (1997). □

## APPENDIX: NOTATION FOR INCREMENTS

Further details concerning increments may be found in Kent and Wood (1997) and Chan and Wood (2000); we only give brief details here.

*The univariate case* $(d=1)$. An increment of order $p \geq 0$ is a finite array of real numbers $\mathbf{a} = \{a_j : -J \leq j \leq J\}$ such that, for all integers $0 \leq r \leq p$,

$$\text{(A.1)} \qquad \sum_{j:\, -J \leq j \leq J} j^r a_j = 0,$$

but

$$\text{(A.2)} \qquad \sum_{j:\, -J \leq j \leq J} j^{p+1} a_j \neq 0.$$

An example of an increment of order $p=0$ is given by

$$\text{(A.3)} \qquad \mathbf{a} = \{a_0, a_1\}, \qquad \text{where } a_0 = -1 \text{ and } a_1 = 1;$$

an example of an increment of order $p=1$ is given by

$$\text{(A.4)} \qquad \mathbf{a} = \{a_{-1}, a_0, a_1\}, \qquad \text{where } a_{-1} = 1, a_0 = -2 \text{ and } a_1 = 1.$$

Note that, here and in the main text, we adopt the convention that any components $a_j$ which are not given explicitly are 0 [so, in (A.3), for example, we have omitted $a_{-1} = 0$].

For an integer $u \geq 1$, the *dilation* $\mathbf{a}^u = \{a_j^u : -Ju \leq j \leq Ju\}$ of an increment $\mathbf{a}$ has elements defined by

$$\text{(A.5)} \qquad a_j^u = \begin{cases} a_{j'}, & \text{if } j = j'u, \\ 0, & \text{otherwise}, \end{cases}$$

where $-Ju \leq j \leq Ju$. It follows immediately from this definition that

$$\sum_j j^r a_j^u = \begin{cases} 0, & \text{if } 0 \leq r \leq p, \\ u^r \sum_j j^r a_j, & \text{if } r \geq p+1, \end{cases}$$

where here and in the main text $\sum_j$ indicates summation over $-Ju \leq j \leq Ju$. Let $\{y_j : j \in \mathbf{Z}\}$ be a sequence of numbers. Then if the increment $\mathbf{a}$ is given by (A.3), it follows that

$$\sum_j a_j^u y_{n+j} = y_{n+u} - y_n,$$

while if $\mathbf{a}$ is given by (A.4), then

$$\sum_j a_j^u y_{n+j} = y_{n+u} + y_{n-u} - 2y_n.$$

*The multivariate case* $(d > 1)$. Let $j = (j[1], \ldots, j[d]) \in \mathbf{Z}^d$ denote a multi-index. We say that $j$ is nonnegative and write $j \geq 0$ if $j[l] \geq 0$ for each $l = 1, \ldots, d$; and if $k$ is another multi-index, we say that $j \leq k$ if $j[l] \leq k[l]$ for each $l$, and write $j < k$ if each inequality is strict. For multi-indices $j$ and $r$, we define

$$|r| = \sum_{l=1}^{d} r[l] \tag{A.6}$$

and

$$j^r = \prod_{l=1}^{d} j[l]^{r[l]}, \tag{A.7}$$

where $0^0 = 1$.

In the $d$-dimensional case, an increment of order $p \geq 0$ is a finite array $\mathbf{a} = \{a_j : -J \leq j \leq J\}$ satisfying (A.1) and (A.2), but with $j$, $J$ and $r$ now multi-indices with $j^r$ defined by (A.7) and each $a_j$ real as before.

The so-called "square" increment in the case $d = 2$, which is an example of an increment of order $p = 1$, is given by

$$\mathbf{a} = \{a_{(0,0)} = a_{(1,1)} = 1, a_{(1,0)} = a_{(0,1)} = -1\}; \tag{A.8}$$

see Chan and Wood (2000).

The dilation $\mathbf{a}^u = \{a_j^u : -Ju \leq j \leq Ju\}$ is defined by (A.5), but with $j$ and $j'$ now multi-indices. It follows immediately from this definition that

$$\sum_j j^r a_j^u = \begin{cases} 0, & \text{if } 0 \leq |r| \leq p, \\ u^{|r|} \sum_j j^r a_j, & \text{if } |r| \geq p+1, \end{cases}$$

where $|r|$ is given by (A.6).

Note that if $\{y_j : j \in \mathbf{Z}^2\}$ is an array of real numbers and the square increment (A.8) is used, then for any $n \in \mathbf{Z}^2$,

$$\sum_j a_j^u y_{n+j} = y_{n+(0,0)} + y_{n+(u,u)} - y_{n+(u,0)} - y_{n+(0,u)}.$$

**Acknowledgments.** We are grateful to the referees for helpful comments which have resulted in an improved presentation.

Department of Statistics
  and Actuarial Science
University of Iowa
Iowa City, Iowa 52242
USA

School of Mathematical Sciences
University of Nottingham
University Park
Nottingham NG7 2RD
United Kingdom
e-mail: andy.wood@maths.nottingham.ac.uk